\documentclass[11pt,leqno]{amsart}
\usepackage[utf8]{inputenc}
\usepackage[all]{nowidow}
\usepackage{amsmath, amssymb, amsthm}
\usepackage{lmodern}
\usepackage{color}
\usepackage{enumerate}
\usepackage{float}
\usepackage{xspace}
\usepackage{longtable}
\usepackage{url}
\usepackage{rotating}
\usepackage[unicode, naturalnames]{hyperref}
\hypersetup{
  colorlinks   = true, 
  urlcolor     = blue, 
  linkcolor    = blue, 
  citecolor   = red 
}

\usepackage[a4paper, inner=1.5in, outer=1in]{geometry}
\usepackage[all]{nowidow}
\usepackage[numbers]{natbib}
\usepackage[symbol]{footmisc}

\newcommand*{\doi}[1]{\href{https://dx.doi.org/#1}{\tt https://doi.org/#1}}

\newcommand{\grp}[1]{\ensuremath{\mathrm{#1}}\xspace}
\newcommand{\aut}[1]{\ensuremath{\mathrm{Aut}(#1)}\xspace}
\newcommand{\inn}[1]{\ensuremath{\mathrm{Inn}(#1)}\xspace}
\newcommand{\symg}[1]{\ensuremath{\mathrm{Sym}(#1)}\xspace}
\newcommand{\autp}[1]{\ensuremath{\mathrm{Aut^+}(#1)}\xspace}
\newcommand{\outp}[1]{\ensuremath{\mathrm{Out^+}(#1)}\xspace}
\newcommand{\sur}[1]{\ensuremath{\mathcal{#1}}\xspace}
\newcommand{\Epi}[1]{\ensuremath{\operatorname{Epi_o}(#1)}\xspace}

\newcommand{\red}[1]{\bgroup\color{red}#1\egroup}
\newcommand{\blue}[1]{\bgroup\color{blue}#1\egroup}

\theoremstyle{plain}
\newtheorem{theorem}{Theorem}[section]
\newtheorem{lemma}[theorem]{Lemma}
\newtheorem{proposition}[theorem]{Proposition}
\newtheorem{corollary}[theorem]{Corollary}
\newtheorem{problem}{Problem}
\theoremstyle{definition}

\newtheorem{example}[theorem]{Example}

\numberwithin{equation}{section}

\begin{document}

\title[Computing equivalence classes\ldots]{Computing equivalence classes of finite\\group actions on orientable surfaces:\\A dynamic survey}

\author[J. Karab\'a\v{s}]{J\'an Karab\'a\v{s}}
\email[J. Karab\'a\v{s}]{jan.karabas@umb.sk}
\author[R. Nedela]{Roman Nedela}
\email[R. Nedela]{nedela@savbb.sk}
\author[M. Skyvov\'a]{M\' aria Skyvov\' a}
\email[M. Skyvov\'a]{mskyvova@ntis.zcu.cz}
\address[J. Karab\'a\v{s}]{Department of Compuer Science, Faculty of Natural Sciences, Matej Bel University, Tajovsk\'eho~40, 97401 Bansk\'a Bystrica, Slovakia}
\address[J. Karab\'a\v{s}, R. Nedela]{Mathematical Institute of Slovak Academy of Sciences, \v{D}umbierska 1, 97411 Bansk\'a Bystrica, Slovakia}
\address[J. Karab\'a\v{s}, R. Nedela, M. Skyvov\'a]{NTIS, Faculty of Applied Sciences, Technick\'a 8, 30601 Plze\v{n}, Czech Republic}
\makeatletter
\@namedef{subjclassname@2020}{%
  \textup{2020} Mathematics Subject Classification}
\makeatother
\date{\today; A dynamic survey, version 2}

\begin{abstract} 
This paper focuses on the classification of
 classes of topological equivalence of finite group actions on Riemann surfaces.
By the Riemann-Hurwitz bound, there are just finitely many groups that act conformally on a closed  orientable surface $\sur{S}_g$ of genus $g\geq 2$. With each such action of a group $\grp{G}$ on $\sur{S}_g$ one can associate the fundamental group $\Gamma=\pi(\sur{O})$ of the quotient orbifold $\sur{O}=\sur{S}_g/\grp{G}$, isomorphic to a Fuchsian group determined completely by orbifold's signature. The Riemann existence theorem reduces the problem of the existence of an action of $\grp{G}$ on $\sur{S}_g$ to a purely group-theoretical problem of deciding whether there is an smooth epimorphism mapping the Fuchsian group $\Gamma$ onto the group $\grp{G}$. Using computer algebra systems such as \textsc{Magma} or GAP, together with the library of small groups, the generation of all finite group actions on a surface of fixed small genus $g\geq 2$ becomes    
almost a routine procedure. The difficult part is to determine the classes of these actions with respect to topological equivalence. To achieve this, one needs to investigate the action of the automorphism group of a Fuchsian group on the set of finite group actions on $\sur{S}_g$ with the corresponding signature.   
In this paper we derive several results on the topological equivalence of finite group actions on Riemann surfaces. As an application, we derive complete lists of finite group actions of genus  $g\leq 9$ distinguished up to the topological equivalence. A summary of the actions can be found in Appendix; the reader interested in more details is referred to the web page 
\url{https://www.savbb.sk/~karabas/science/discactions.html}
It is expected that we will be able to extend the list to higher genera, refreshed partial
results are available on the web page. The following text is an extended
version of the paper \cite{KNS24}.
\end{abstract}
\maketitle
\thispagestyle{empty}

\section{Introduction}

\noindent{}The investigation of surface symmetries is a long-standing project started in ancient Greece. By the beginning of 20th century
the classification of the spherical and wallpaper groups, which are the groups of isometries of the sphere and of the Euclidean plane, had been established. Groups acting on the projective plane, the torus, and the Klein bottle are understood as well, since these are quotients of the spherical and wallpaper groups. In contemporary mathematics,
the symmetries of surfaces of higher genera are studied within the frame of the theory of Riemann surfaces and of algebraic curves. With each orientation-preserving action of an group $\grp{G}$ on a orientable surface $\sur{S}_g$ of genus $g$, one can associate the \emph{signature} 
$(\gamma;m_1,m_2,\dots,m_r)$, where $\gamma$ is the genus of the quotient surface, $r$ is the number of singular orbits, and $m_i$, $i=1,\dots,r$, are the orders of the point stabilizers in the corresponding singular orbits (the point stabilizers are necessarily cyclic). 
The signature, genus $g$, and order of the group $\grp{G}$ are related by the Riemann-Hurwitz equation~\eqref{eq:rh}. As a consequence, the order of $\grp{G}$ is related via the Hurwitz bound, $|\grp{G}|\leq 84(g-1)$ for $g>1$, see e.g.~\cite[Theorem~5.11.1]{js87}. It follows that for genera greater than one, there are only
finitely many finite groups acting on these surfaces. Thus, determining all of them is a finite problem. On the other hand, groups that are small with respect to Hurwitz bound
have usually a huge number of actions even on surfaces
of small genera. Hence it makes sense to consider proper equivalence relations and classify the actions up to the chosen
equivalence. An excellent account on equivalence relations of group actions on surfaces, together with many open problems, can be found in recent papers \cite{broughton2022equivalence, broughton2022future}. 

One of the most important relations defined on actions of groups on orientable surfaces is the topological
equivalence defined as follows.
For an orientable surface $\sur{S}$ denote by $\operatorname{Hom}^+(\sur{S})$ the group of orientation-preserving homeomorphisms. Two actions of a group $\grp{G}$, given by embeddings 
$\varepsilon_i\colon \grp{G}\to \operatorname{Hom}^+(\sur{S}_i)$, $i=1,2$, on possibly different surfaces
$\sur{S}_1$, $\sur{S}_2$ of
the same genus, are \emph{topologically equivalent} if there is an intertwining, orientation
preserving homeomorphism $h\colon \sur{S}_1\to \sur{S}_2$ and an automorphism $a\in  \aut{\grp{G}}$ such that
$\varepsilon_2(g)=h\varepsilon_1(a(g))h^{-1},$ for every $g\in \grp{G}$. As it is noted in \cite[p.23]{broughton2022equivalence}, ``topological equivalence of actions is important since the action classes are directly linked to conjugacy classes of finite subgroups of the mapping class group, in a 1-1 fashion, and to the strata of the branch locus of moduli space in an almost 1-1 fashion''.

\begin{problem} \label{Problem} For a fixed integer $g>1$, derive the list of actions of
finite groups  on an orientable surface of genus $g$, distinguished up to topological equivalence.
\end{problem}

A complete solution of Problem~\ref{Problem} is known only up to genus $4$,  see Broughton~\cite{broughton1991classifying} for genera $2$ and $3$, and Bogopolsky~\cite{bogopolski1997classification} for genus four.
Incomplete lists for genera $4$ and $5$ were derived by  Kimura~\cite{Kimurag4}, and by Kuribayashi and Kimura~\cite{Kuribayashi}.  To determine the classes of topological equivalence of finite groups of genus $g$ seems to be a hard problem, both from theoretical and practical point of view. Broughton in \cite{broughton2022equivalence} discusses several related
equivalences of group actions which may be used to approximate the topological equivalence. In particular,  Breuer in \cite{breuer2000characters} explains
how one can employ computations with group characters to determine the equivalence
relation $\sim_{\mathcal{H}^1}$ (following the notation from   
\cite{broughton2022equivalence}). The numbers of $\sim_{\mathcal{H}^1}$-classes
are computed up to genus $48$, see \cite[p.95]{breuer2000characters}.
By Broughton 
\cite[p.28]{broughton2022equivalence}, two actions that are topologically equivalent  
are $\sim_{\mathcal{H}^1}$-equivalent as well. Unfortunately, in many instances of Problem~\ref{Problem} the topological equivalence differs significantly from $\sim_{\mathcal{H}^1}$-equivalence, see discussion in Section~\ref{sec:class} for details. Paulhus in \cite{paulhus2019dbpaper} introduces a database of group actions, however, the actions are not distinguished up to topological equivalence. Moreover, we have checked that the database \cite{paulhusdb}  is not complete, see Section~\ref{sec:class} for some details.

An action of a finite group $\grp{G}$ on (an orientable) surface of genus $g$ determines
the quotient \emph{orbifold} $\sur{S}_g/\grp{G}$ with the fundamental group $\Gamma$. On the other hand, $\Gamma$ acts as a group of orientation-preserving symmetries of the universal cover $\sur{U}$ of the surface $\sur{S}_g$. In this setting, $\Gamma$ is known under the name \emph{Fuchsian group}.  
The surface $\sur{S}_g$ and the group $\grp{G}$ can be obtained as quotients $\sur{U}/K$ and $\Gamma/K$ by some normal subgroup $K\lhd\Gamma$, respectively. The group $K$ is a surface group isomorphic to $\pi_1(\sur{S}_g)$.  
The Fuchsian group
$\Gamma$ is completely determined by the signature
$(\gamma;m_1,\dots,m_r)$, therefore, we use the notation
$\Gamma=\grp{F}(\gamma;m_1,\dots,m_r)$.
  By the Riemann existence theorem (Theorem~\ref{thm:existence}), an action of an abstract finite group $\grp{G}$ on $\sur S_g$ with a prescribed signature exists if and only if there is a smooth
epimorphism $\eta\colon\Gamma\to \grp{G}$. An epimorphism is \emph{smooth}
if for each element $x$ of finite order $|\eta(x)|=|x|$. Equivalently, smooth epimorphisms are called \emph{surface kernel} epimorphisms \cite{Lloyd}. 
In view of the above, the Riemann existence theorem yields that the classification of discrete groups of a given genus can be seen as a purely group-theoretical problem. More precisely, Lloyd proved in~\cite{Lloyd} that
two actions of a group $\grp{G}$ with given signature, given by associated epimorphisms $\eta_1,\eta_2\colon \Gamma\to\grp{G}$, are topologically equivalent if and only if there
exist $a\in\aut{\grp G}$ and orientation-preserving automorphism $\alpha\in \autp{\Gamma}$ such that 
$\eta_2=a\eta_1\alpha$. 
It follows that in order
to determine the equivalence classes it is necessary 
to study the action of the \emph{group of orientation-preserving automorphisms} $\autp{\Gamma}$ of the Fuchsian group $\Gamma$
on the set of all smooth epimorphisms $\Gamma\to \grp{G}$. A particular subproblem is to find a proper finite set of generators of $\autp{\Gamma}$. The problem of determining a finite generating sets of $\autp{\Gamma}$ is related to the problem of determining the generating sets of mapping class groups of two-dimensional orbifolds, a topic intensely studied within the framework of geometric topology. In the cases we are interested in, the mapping class group is isomorphic to the outer automorphism group of $\Gamma$. Therefore, finding a finite generating set for $\autp{\Gamma}$ and finding the generating set for the corresponding mapping class group are equivalent problems; see \cite{bridson2005,broughton2022future} for details. Various generating sets of mapping class groups of surfaces are determined 
in \cite{lickorish1964finite, mccool1975some, humphries1979generators, hatcher1980presentation, wajnryb1983simple, korkmaz2005generating, stukow2014finite}. In \cite{gervais2001finite, labruere2001presentations} one can find a finite set of generators for mapping class groups of punctured surfaces.

 A significant part of the paper deals with symmetries of Fuchsian groups with planar signatures. Following Zieschang et al.~\cite{zieschang2006surfaces},
a Fuchsian group is \emph{planar} if its signature is of the form $(0;m_1,m_2,\dots,m_r)$.  Finite sets of generators and presentations of $\autp{\Gamma}$, for planar $\Gamma$ and $r\leq 6$ can be found in \cite{tap88}. In Section~2
we introduce basic concepts and relevant general results.  In Section~3 we collect and organise known results on the automorphism group of a planar Fuchsian group derived by Zieschang \cite{zieschang2006surfaces}, Birman \cite{birman1974braids}, and Chow \cite{chow2002algebraical}. As a consequence, in Section~4 we are able to formulate and prove Theorem~\ref{thm:equiactions} characterising the action of $\autp{\Gamma}$ on the \emph{sets of smooth epimorphisms}~$\Epi{\Gamma,\grp{G}}$ for Fuchsian groups $\Gamma$ with planar signatures. Based on the theoretical part of the paper we were able to implement an algorithm determining representatives of equivalence classes of finite groups with planar signatures  for a fixed genus. 
Further, employing the generating set of the automorphism group of a surface group derived by McCool in \cite{mccool1975some} we  implemented the algorithm computing equivalence
classes of actions with signatures of type $(\gamma;-)$, $\gamma\geq 1$.
The two algorithms are complemented by an approximation algorithm producing  an equivalence relation refining the topological equivalence for actions with nonplanar signature. The approximation uses automorphisms of Fuchsian groups introduced in~\cite{harvey1971}. Actual outputs of programs and other related information can be found online~\cite{karaweb}. In particular, in this paper we deal with actions of finite groups of genera $2\leq g\leq 9$, see Section~\ref{sec:class} and Appendix for more details.

\section{Finite group actions on surfaces}
\label{sec:Gactions}

\noindent{}The following set of transformations, known as the \emph{group of M\"obius transformations}, 
maps the upper half-plane $\sur{U}=\{ z \in \mathbb{C} \ |\ \operatorname{Im}(z) > 0 \}$  onto itself by taking  $z\mapsto w(z)=\frac{az+b}{cz+d}$ for $a,b,c,d \in \mathbb{R}$, $ad-bc=1$. The group, known as $\grp{LF}(2,\mathbb{R})$ or $\grp{PSL}(2,\mathbb{R})$, can be topologised by
taking the usual metric topology in $\mathbb{R}^4$ and identifying $w$ with the equivalence
class of $(a, b, c, d)\in \mathbb{R}^4$, $(a, b, c, d) \sim \lambda(a, b, c, d)$ for all $\lambda \in \mathbb{R}\setminus \{ 0\}.$
A discrete subgroup of $\grp{PSL}(2, \mathbb{R})$ is called a Fuchsian group.
Fuchsian groups for which the orbit space  is compact will be of
interest here, so henceforth `Fuchsian group' will mean `Fuchsian group
with compact orbit space'. 

A torsion-free Fuchsian group is a \textit{surface group}. Surface groups are isomorphic to fundamental groups of compact connected orientable surfaces. It is known that a Fuchsian group $\Gamma$ always contains a normal surface subgroup $\grp{K}$ of finite index \cite{nielsen1927}. On the other hand, every normal torsion-free subgroup of finite index in $\Gamma$ is a surface group.  Recall that a discrete group of automorphisms of a connected compact closed surface is finite. By Zieschang~\cite{zieschang1966automorphismen}, every such group is a quotient of a subgroup of the group of M\"obius transformations of the upper half-plane \sur{U}. 

Let $\sur{S}_g$ be an orientable surface of the genus $g$, and let $\operatorname{Hom}^+(\sur{S}_g)$ be its group of orientation-preserving
homeomorphisms. We say that the group $\grp{G}$ acts on the surface $\sur{S}_g$ if there is a monomorphism $\varepsilon\colon \grp{G}\to \operatorname{Hom}^+(\sur{S}_g)$.
Every action of \grp{G} on $\sur{S}_g$ is determined by the pair of Fuchsian groups $\grp{K}\lhd \Gamma<\grp{PSL}(2,\mathbb{R})$
that act discontinuously on $\sur{U}$ and by an epimorphism $\eta\colon \Gamma\to \grp{G}$ with kernel $\grp{K}$, where $\grp{K}$ 
is a surface group. Such an epimorphism is 
smooth by definition.
The epimorphism $\eta$ is constructed from $\varepsilon$ and from a homeomorphism of $\sur{U}/\grp{K}\cong \sur{S}_g$.
More precisely, let \grp{G} be a finite group. The smooth epimorphism $\eta\colon \grp{\Gamma} \to \grp{G}$ determines a discrete action of the group \grp{G} on the quotient Riemann surface $\sur{U}/\ker\eta$. 
 For every finite group \grp{G} acting as a group of automorphisms of a compact orientable surface~$\sur{S}$ there is a Fuchsian group $\Gamma$ and an smooth epimorphism $\eta\colon \Gamma \to \grp{G}$, such that the action $(\grp{G},\sur{S})$ is equivalent to the action $(\Gamma/\ker\eta,\ \sur{U}/\ker\eta)$.

The group $\Gamma = \grp{F}(\gamma;m_1, m_2,\ldots, m_r)$ has the (canonical) presentation

\begin{align}\label{eq:fuchs}
\langle x_1,x_2,\ldots,x_r,a_1,\ldots,a_\gamma,b_1,\ldots,b_\gamma\ |\ &x_1^{m_1} = \cdots = x_r^{m_r}=1,\ \prod_{i=1}^{\gamma}[a_i,b_i]\prod_{j=1}^{r}x_j = 1\rangle,
\end{align}
where $1<m_1\leq m_2\leq\dots m_r$ are integers, $0\leq\gamma\leq g$ is the genus of the quotient
surface $\sur{S}_g/\eta(\Gamma)$, and $m_i$ are the branch indices of the (regular) branched covering 
$\sur{S}_g\to \sur{S}_g/\eta(\Gamma)$, such that $m_i$ divides $|\grp{G}|$ for $i=1,2,\dots,r$.
All the integer parameters in the presentation~\eqref{eq:fuchs} are related with the genus $g$ of the surface $\sur{S}_g$  by the Riemann-Hurwitz equation

\begin{align}\label{eq:rh}
2-2g {}= {}&|\grp{G}|\left(2-2\gamma-\sum^r_{i=1}\left(1-\frac{1}{m_i}\right)\right);\\& m_i\geq 2,\quad m_i\mid|\grp{G}|,\quad i=1,\dots,r.\notag
\end{align}

As a consequence, one finds that if $g\geq 2$, then the order of the group $\grp{G}$ is bounded by the Hurwitz bound, $|\grp{G}|\leq 84(g-1)$. It follows that there are finitely many groups acting on the surface $\sur{S}_g$ of genus $g\geq 2$. The integers $\gamma$, $m_1,\dots,m_r$ determined by the action of $\grp{G}$, usually written as vector $(\gamma;m_1,\dots, m_r)$, will be called the \emph{signature} or the \emph{branch data}.
Consider the inverse problem:
Given the genus $g\geq 2$, let $|\grp{G}|$ and a signature $(\gamma;m_1,\dots, m_r)$ satisfy the Riemann-Hurwitz equation.
Is there a group $\grp{G}$ of order $|\grp{G}|$ that acts on the orientable surface $\sur{S}_g$ of the genus $g$ such that the Riemann-Hurwitz equation is satisfied for the signature $(\gamma;m_1, m_2,\ldots,m_r)$?
It turns out that this is not always true. A complete answer is given by the following theorem.

\begin{theorem}[Riemann existence theorem~\cite{broughton1991classifying, harvey1966cyclic}]\label{thm:existence} Group $\grp{G}$ acts on a surface $\sur{S}_g$ of genus $g$ with branch data 
$(\gamma;m_1,\dots,m_r)$ if and only if there is an smooth epimorphism 
$\grp{F}(\gamma;m_1,\dots,m_r)\to \grp{G}$.
\end{theorem}

Recall that the epimorphism $\eta\colon\Gamma\to\grp{G}$ is smooth if and only if $|\eta(x_i)|=m_i=|x_i|$, for $i=1,2,\dots,r$. We say that a pair $(\Gamma,\grp{G})$,
where $\grp{G}$ is a finite group and $\Gamma=\grp{F}(\gamma;m_1,\dots,m_r)$ is a Fuchsian group, is \emph{numerically $g$-admissible} if the corresponding parameters satisfy (\ref{eq:rh}).
If $\Epi{\Gamma,\grp{G}}$ is non-empty, the we say
that $(\Gamma,\grp{G})$ is \emph{$g$-admissible}.


The following theorem gives an algebraic criterion for the equivalence
of actions in terms of the associated smooth epimorphisms.
\begin{theorem}[{\cite[Theorem~4]{Lloyd}}]
\label{Lloyd}
 Let $\grp{\Gamma}$ be a Fuchsian group. Given a finite group \grp{G}, let $\eta_1, \eta_2 \in \Epi{\grp{\Gamma},\grp{G}}$ be smooth epimorphisms, and let $(\grp{G},\sur{S})$ and $(\grp{G},\sur{S}')$ be
the corresponding actions. Then $(\grp{G},\sur{S})$ and
$(\grp{G},\sur{S}')$ are topologically equivalent if and only if there exists $\alpha \in \autp{\grp{\Gamma}}$ and $a \in \aut{\grp{G}}$
such that $\eta_2=a\eta_1 \alpha$.
\end{theorem}
The smooth epimorphisms $\eta_1$, $\eta_2$ satisfying the conditions in Theorem~\ref{Lloyd} will be called  \emph{equivalent}, in short, we write $\eta_1 \sim \eta_2$. If $\eta_2=a\eta_1$ for some $a\in\aut{\grp{G}}$, then $\eta_1\cong \eta_2$ are \emph{$\aut{\grp{G}}$-equivalent}.

Let $\Gamma$ be a Fuchsian group given by presentation (\ref{eq:fuchs}) and $\grp{G}$ be a finite group.
Every epimorphism $\eta\in\Epi{\Gamma,\grp{G}}$ is determined by the images
of the $r+2\gamma$ generators from presentation (\ref{eq:fuchs}) written as a vector
$$ \vec\eta=(\bar x_1,\dots,\bar x_r,\bar a_1,\dots,\bar a_\gamma,\bar b_1,\dots,\bar b_\gamma),$$
where $\bar x_i=\eta(x_i)$, $\bar a_j=\eta(a_j)$, $\bar b_j=\eta(b_j)$, for $i=1,\dots,r$
and $j=1,\dots,\gamma$. Hence the set $\Epi{\Gamma,\grp{G}}$ of smooth epimorphisms $\Gamma\to \grp{G}$ can be identified with the set of above vectors of length $r+2\gamma$.
Let
\[GEN=(x_1,\dots, x_r, a_1,\dots, a_\gamma, b_1,\dots, b_\gamma)\]
be a vector (treated as an array) of generators of a Fuchsian group $\Gamma$. Then we can express
a modified version of Theorem~\ref{Lloyd} as follows.

\begin{corollary}[{\cite[p.22]{broughton2022equivalence}}]
The classes of topological equivalence are in one-to-one correspondence
with the orbits  of $\aut{\grp{G}}\times \autp{\Gamma}$ on $\Epi{\Gamma,\grp{G}}$
given by $(a,\alpha)\cdot\vec\eta_1=\vec\eta_2$, where the $k$-th entry of $\vec\eta_2$ is $\vec\eta_2[k]=a(\eta_1(\alpha(GEN[k])))$, for $k=1,\dots,r+2\gamma$.
\end{corollary}
  
Let $\eta\in \Epi{\Gamma,\grp{G}}$.
Write the multi-set $\{\eta(x_i)|\ i=1,\dots,r\}$
in the form $\{y_1^{n_1},\dots,y_k^{n_k}\}$, where
$n_1\leq n_2\leq\dots\leq n_k$, where every $n_i$ is the multiplicity of an element $y\in \{\eta(x_i)|\ i=1,\dots,r\}$. Clearly, $r=n_1+n_2+\dots+n_k$.
Denote 
$s(\eta)=(n_1,\dots,n_k)$. We show that $s(\eta)$
is an invariant of a class of topological equivalence
provided $G$ is abelian.

\begin{lemma}\label{lem:seta} Let $(\Gamma,G)$ be a $g$-admissible
pair and let $G$ be abelian. If $\eta_1,\eta_2\in
\Epi{\Gamma,\grp{G}}$ are topologically equivalent, then $s(\eta_1)=s(\eta_2)$.
\end{lemma}
\begin{proof}
By Theorem~\ref{Lloyd} there exist 
$\alpha\in\autp{\Gamma}$ and $a\in \aut{\grp{G}}$
such that $\eta_2=a\eta_1\alpha$. By \cite[Theorem~5.8.6]{js87}
there exists a permutation $\mu\in \symg{r}$ and $g_i\in\Gamma$ such that $\alpha(x_i)= g_ix_{\mu(i)}g_i^{-1}$, for  $i=1,\dots,r$. It follows that
\[\eta_2(x_i)=a(\eta_1(\alpha(x_i)))=
a(\eta_1(g_ix_{\mu(i)}g_i^{-1}))=a(\eta_1(g_i)\eta_1(x_{\mu(i)}){\eta_1(g_i)}^{-1}=a(\eta_1(x_{\mu(i)})).\]
Since  
$\{\eta_2(x_1),\dots,\eta_2(x_r)\}=\{a\eta_1(x_{\mu(1)}),\dots,a\eta_1(x_{\mu(r)})\}$, we have $s(\eta_1)=s(\eta_2)$.
\end{proof}

\begin{lemma}\label{lem:normaliser} Let $(\Gamma,G)$ be a $g$-admissible
pair. If $\eta_1,\eta_2\in
\Epi{\Gamma,\grp{G}}$ are topologically equivalent, then there exists $a\in \aut{\grp{G}}$ such that $N_\grp{G}(\langle \eta_1(x_1),\dots,\eta_1(x_r)\rangle)=a N_\grp{G}(\langle \eta_2(x_1),\dots,\eta_2(x_r)\rangle)$, where $N_\grp{G}(\grp{H})$ denotes the
normaliser of $\grp{H}\leq \grp{G}$ in $\grp{G}$.
\end{lemma}
\begin{proof}
By Theorem~\ref{Lloyd} there exist 
$\alpha\in\autp{\Gamma}$ and $a\in \aut{\grp{G}}$
such that $\eta_2=a\eta_1\alpha$. Set 
$K=\langle x_1,\dots,x_r\rangle\leq\Gamma$.
By \cite[Theorem~5.8.6]{js87} we have 
$\alpha(N_\Gamma(K))=N_\Gamma(K)$.

It follows that
$$\eta_2(N_\Gamma(K))=a(\eta_1(\alpha(N_\Gamma(K)))=
a(\eta_1(N_\Gamma(K)))= a N_\grp{G}(\langle \eta_1(x_1),\dots,\eta_1(x_r)\rangle).
$$
On the other hand, $\eta_2(N_\Gamma(K))=N_\grp{G}(\langle \eta_2(x_1),\dots,\eta_2(x_r)\rangle)$,
and we are done.
\end{proof}

\section{Automorphism group of a planar Fuchsian group}
\label{sec:planar}
\noindent{}This section is aimed to derive a finite generating set of $\autp{\Gamma}$ for a \emph{planar Fuchsian group}. 
By a planar Fuchsian group $\Gamma = \grp{F}(0;m_1,m_2,\ldots,m_r)$ we mean the group with presentation
\begin{align}\label{eq:planarfuchsian}
\langle x_1,x_2,\ldots,x_r\ |\ x_1^{m_1} = \cdots = x_r^{m_r}=1, 
\prod_{j=1}^{r}x_j = 1\rangle,\ \ 2+\sum_{i=1}^r\frac{1}{m_i}<r.
\end{align}
In what follows, in Sections~\ref{sec:planar} and~\ref{sec:mappings}, the symbol $\Gamma$ will denote a planar Fuchsian group.
 
Denote by $\symg{r}$ the symmetric group that acts on the set $\{1,2,\dots,r\}$ and by
 $\mathcal{F}_r=\mathcal{F}(\tilde x_1, \dots, \tilde x_r)$ the free group of rank $r$ generated by the free generators 
$\tilde x_1, \dots, \tilde x_r$.
Let $\grp{B}_r=\langle\sigma_1, \dots,\sigma_{r-1} \rangle$, be the \emph{braid group} with the defining relations
$\sigma_i\sigma_{i+1}\sigma_i=\sigma_{i+1}\sigma_{i}\sigma_{i+1}$, for $i=1,\dots,r-1$, and $\sigma_i\sigma_j=\sigma_j\sigma_i$, for $|i-j|\geq 2$. 
Let $\nu\colon \grp{B}_r \to \symg{r}$ be the homomorphism defined by 
$\nu(\sigma_i)=(i,i+1)$. Since $\grp{B}_r$ acts on $\symg{r}$ through a transitive homomorphism $\nu$ and for $r\geq 3$ the image contains
a $3$-cycle, the image $\nu(\grp{B}_r)$ is $\symg{r}$.
Hence, $\nu$ is a group epimorphism $\grp{B}_r\to \symg{r}$.
The kernel $\ker(\nu)$ is called the \emph{pure braid group} 
$\grp{P}_r$. It is known that $\grp{P}_r$ is generated by the elements 
$$A_{i,j}={\sigma^{-1}_{j-1}}{\sigma^{-1}_{j-2}}\dots\sigma^{-1}_{i+1}{\sigma_i}^2\sigma_{i+1}\dots\sigma_{j-1},$$ 
for $1\leq i<j\leq r$, see~\cite{chow2002algebraical}. 

Denote by $\delta \colon \grp{B}_r \hookrightarrow \aut{\mathcal{F}_r}$ the embedding determined by $\sigma_i\mapsto c_i$, where \[c_i:\tilde x_i\mapsto \tilde x_i\tilde x_{i+1}\tilde x_i^{-1},\ \tilde x_{i+1}\mapsto \tilde x_i\text{ and }\tilde x_j\mapsto \tilde x_j
\text{ for } 
j\notin\{i,i+1\}.\]
We now consider the following natural question: How do we recognize the image $\delta(\grp{B}_r)$ in $\aut{\mathcal{F}_r}$?
This was answered by Birman 1974 \cite{birman1974braids} and Chow 1948 \cite{chow2002algebraical} as follows.

\begin{theorem}[\cite{chow2002algebraical}]\label{thm:braid} An automorphism 
$\beta\in\delta(\grp{B}_r)$ if and only if
there exists a permutation $\mu\in \symg{r}$ and elements $\lambda_i\in{\mathcal F}_r$ such
that 
\begin{itemize}
\item[(i)] $\beta(\tilde x_i)=\lambda_i \tilde x_{\mu(i)}\lambda^{-1}_i$, 
$i=1,2,\dots,r$, and
\item[(ii)] $\beta(\tilde x_1\tilde x_2\dots\tilde x_r)=
\tilde x_1\tilde x_2\dots\tilde x_r$.
\end{itemize}
\end{theorem}

By \cite[Theorem~3]{zieschang1966automorphismen}, every automorphism of a Fuchsian group $\Gamma$ is induced by an automorphism of
the free group $\mathcal{F}_r$.
The following theorem characterizing the  
 lift  
$\tilde A(\Gamma)\leq \aut{\mathcal{F}_r}$ of 
$\autp{\Gamma}$ is a restricted version of \cite[Theorem~3]{zieschang1966automorphismen}. 

\begin{theorem}[{\cite[Theorem~8]{zieschang1966automorphismen}}]\label{thm:lift} If an automorphism $\tilde\phi\in\tilde A(\Gamma)\leq\aut{\mathcal{F}_r}$, then there exists a permutation $\mu\in \symg{r}$ satisfying $|x_{\mu(i)}|=m_i=|x_i|$ and elements $\lambda_i$, $\lambda \in \mathcal{F}_r$, for $i=1,\dots,r$,  such
that
\begin{itemize} 
\item[(i)] $\tilde\phi(\tilde x_i)=\lambda_i \tilde x_{\mu(i)}\lambda^{-1}_i$, 
$i=1,2,\dots,r$, and
\item[(ii)] $\tilde\phi(\tilde x_1\tilde x_2\dots\tilde x_r)=
\lambda\tilde x_1\tilde x_2\dots\tilde x_r\lambda^{-1}$,
\item[(iii)] $\omega(\lambda)=1$, where $\omega\colon \mathcal{F}_r\to\Gamma$ is the natural projection. 
\end{itemize}
Furthermore, given $\tilde\phi\in\tilde A(\Gamma)$ determined by $\mu\in \symg{r}$, $\lambda$ and $\lambda_i$ ($i=1,\dots,r$)
projects to a unique $\phi\in\autp{\Gamma}$. 
\end{theorem}

The following lemma is proved in \cite{tap88}.

\begin{lemma}[{\cite[Lemma~2.1]{tap88}\label{lem:inn}}] The lift $\tilde I(\Gamma)\leq \aut{{\mathcal F}_r}$ of the group
of inner automorphisms $\inn{\Gamma}$ is a subgroup of $\delta(P_r)$.
\end{lemma}

Using the aforementioned statements, we can prove the following proposition.

\begin{proposition}\label{prop:nustar}Let $\Gamma$ be a planar Fuchsian group. Then $\autp{\Gamma}$ acts on the set $\{1,2,\dots,r\}$ by a homomorphism $\nu^*\colon \autp{\Gamma}\to \symg{r}$ defined by setting
$\nu^*(\phi)=\nu(\delta^{-1}(\tilde\phi))$, where $\tilde\phi$ is a lift of $\phi\in \autp{\Gamma}$.
\end{proposition}

\begin{proof} It is sufficient to prove that $\nu^*$ is well defined.  According to Lemma~\ref{lem:inn}, the group $\tilde I(\Gamma)$ of lifts of inner automorphisms of $\Gamma$ is a subgroup of $\delta(P_r)=\delta(\ker(\nu))$.
By Theorem~\ref{thm:lift} and Lemma~\ref{lem:inn} we have $\tilde A(\Gamma)\leq \tilde I(\Gamma)\cdot \delta(B_r)=\delta(B_r)$. It follows that two lifts  $\tilde \phi_j$, $j=1,2$, of $\phi\in\autp{\Gamma}$  are described by the elements $\mu_j\in\symg{r}$
and $\lambda_{j,i}\in \mathcal{F}_r$, $i=1,\dots,r$. Since $\tilde\phi_j$ are lifts of the same automorphism $\phi$,
the product $\kappa=\tilde\phi_1\cdot\tilde\phi_2^{-1}$ is a lift of identity. Hence $\mu_1=\mu_2$ and 
$\kappa\in \delta(P_r)$. It follows that $\tilde\phi_1=\kappa\tilde\phi_2$. Since $\delta^{-1}(\kappa)\in\ker(\nu)$, we have $\nu\delta^{-1}(\tilde\phi_1)=\nu\delta^{-1}(\tilde\phi_2)$. It follows that $\nu^*(\phi)$ does not
depend on the choice of the lift $\tilde\phi$.
\end{proof}

Let $\mathcal P$ be the partition
of the index set $\{1,2,\dots,r\}$ given by the equivalence
$i\cong j$ iff $m_i=m_j$. We assume that the integers $m_i$ are ordered
in non-decreasing order. Assume that there are $\ell\geq 1$ classes
of $\mathcal P$. Then there exist integers $0=r_0 < r_1<r_2<\dots <r_\ell=r$ such that
each part of $\mathcal P$ has form
$${\mathcal P}_j=\{ r_{j-1}+1,\dots,r_j\},$$
where 
$r_j-r_{j-1}=|{\mathcal P}_j|$, for some $j$, $\ell\geq j\geq 1$. We say that a permutation $\mu\in\symg{r}$ 
is \emph{$\mathcal{P}$-invariant} if $\mu(\mathcal{P}_j)=\mathcal{P}_j$, for $j=1,\dots,\ell$.
The set of all $\mathcal{P}$-invariant permutations forms a subgroup $\symg{r,{\mathcal P}}\leq\symg{r}$
isomorphic to the direct product $\prod_{j=1}^\ell \symg{r_j-r_{j-1}}$.

\begin{theorem}\label{thm:kerim} Let $\Gamma$ be a planar Fuchsian group. Let $\omega\colon \tilde A(\Gamma)\to \autp{\Gamma}$ be the natural projection. Then for the epimorphism \mbox{$\nu^*\colon \autp{\Gamma}\to \symg{r}$} 
the following statements hold:
\begin{itemize}
\item[(i)] The kernel $\ker{\nu^*}$ is the image $P^*_r=\omega\delta(P_r)$ of the pure braid group;
\item[(ii)] The image $\nu^*(\autp{\Gamma})=\symg{r,{\mathcal P}}\cong \prod_{j=1}^\ell \symg{r_j-r_{j-1}}$
is the group of $\mathcal P$-invariant permutations;
\item[(iii)] $\autp{\Gamma}=\omega\delta\nu^{-1}(\symg{r,{\mathcal P}})$.
\end{itemize}
\end{theorem}

\begin{proof} By  definition, a generator $c_i$, $i=1,2,\dots,r-1$,
of $\delta(\grp{B}_r)$ is in $\tilde A(\Gamma)$ exactly
when $i\neq r_j$ for some $j$. Therefore, $\symg{r,{\mathcal P}}=\nu^*(\autp{\Gamma)}$.

Let $\phi\in\ker(\nu^*)$. Then $\nu(\delta^{-1}(\tilde \phi))=1$.
It follows that $\tilde\phi\in\delta(P_r)$, and consequently, $\tilde\phi\in \tilde I(\Gamma)\cdot \delta(P_r)$, where $\tilde{I}(\Gamma)$ is the lift of the inner automorphisms (see Lemma~\ref{lem:inn}).
By \cite[Lemma~2.1]{tap88} we have $\tilde I(\Gamma)\cdot \delta(P_r)=\delta(P_r)$, and we have $\ker(\nu^*)\leq \omega\delta(P_r)=P^*_r$. Since $P_r^*\leq\ker(\nu^*)$ holds trivially, $P_r^*=\ker(\nu^*)$.

 By Theorem~\ref{thm:lift} and by  \cite[Lemma~2.1]{tap88} we have $\tilde A(\Gamma)=\tilde I(\Gamma)\cdot \delta\nu^{-1}(\symg{r,{\mathcal P}})= \delta\nu^{-1}(\symg{r,{\mathcal P}})$. It follows
 that $\autp{\Gamma}=\omega \delta\nu^{-1}(\symg{r,{\mathcal P}})$.
\end{proof}

\begin{corollary} For a permutation $\mu\in\symg{r,{\mathcal P}}$ choose an element $x_\mu\in\autp{\Gamma}$ with the property $\nu^*(x_\mu)=\mu$. Then $\autp{\Gamma}=\cup_\mu x_\mu P_r^*$. In particular, the index
$[\autp{\Gamma}:P_r^*]=\prod_{j=1}^\ell (r_j-r_{j-1})!$.
\end{corollary}

An explicit set of generators of the group $P^*_r\leq \autp{\Gamma}$ can be derived
from the generating set $\{A_{s,t}\}$ of the pure braid group; see \cite{tap88}.
For $1 \leq s < t \leq r$ we set
\begin{align}\label{eq:purebraid}
\mathcal{A}_{s,t}(x_i) &=\left\{\begin{array}{llll}
x_i, & &  i<s \text{ or } t<i,\\
\rule{0pt}{16pt}x_i^{(x_sx_t)}&=x_sx_tx_ix_t^{-1}x_s^{-1}, & s=i, \\
\rule{0pt}{16pt}x_i^{x_s} &= x_sx_ix_s^{-1}, & t=i, \\
\rule{0pt}{16pt}x_i^{[x_s,x_t]}&=x_sx_tx_s^{-1}x_t^{-1}x_ix_tx_sx_t^{-1}x_s^{-1}, & s<i<t  
\end{array}\right. 
\end{align}

\begin{corollary} Let $\omega\colon \mathcal{F}_r\to {\Gamma}$ be the natural projection
$\tilde x_i\to x_i$, $i=1,\dots,r$. Then
the group $\autp{\Gamma}$ is generated by the set 
$$\{{\mathcal A}_{s,t}:\ 1\leq s<t\leq r\}\cup \{\omega c_i: 
i\in \{1,2,\dots,r\}\setminus\{r_1,\dots,r_\ell\}\}.$$
\end{corollary}
In what follows, we set $\mathcal{H}_i = \omega c_i\in\autp{\Gamma}$, for $i\in \{1,2,\dots,r\}\setminus\{r_1,\dots,r_\ell\}$.

\section{Action of $\autp{\Gamma}$  on $\Epi{\Gamma,\grp{G}}$}
\label{sec:mappings}

\noindent{}Let $\Gamma=\grp{F}(0;m_1,m_2,\ldots,m_r)$ be a planar Fuchsian group, and let $\grp{G}$ be a finite group.
The action of 
$\autp{\Gamma}$ on $\Epi{\Gamma,\grp{G}}$ is given by
$\eta\mapsto \eta\alpha$, for $\eta\in \Epi{\Gamma,\grp{G}}$ and
$\alpha\in\autp{\Gamma}$. In what follows we use
the representation of epimorphisms from $\Epi{\Gamma,\grp{G}}$  by the vectors of images of the generators of $\Gamma$ introduced in Section~\ref{sec:Gactions}.

We call  \emph{ vertical action}  the action of  the image $P^*_r\leq \autp{\Gamma}$
of the pure braid group on $\Epi{\Gamma,\grp{G}}$ . It is generated by the permutations $\alpha_{s,t}$ corresponding to the  generators ${\mathcal A}_{s,t}$:
\begin{equation}\alpha_{s,t}\colon(\eta(x_1),\dots,\eta(x_r))\mapsto ({\eta(\mathcal A}_{s,t}( x_1)),\dots,\eta({\mathcal A}_{s,t}(x_r))),\text{ for } 1\leq s<t\leq r.\end{equation}

The images ${\mathcal A}_{s,t}(x_i)$ are defined by \eqref{eq:purebraid}. It follows
from Theorem~\ref{thm:kerim}(iii) that
the action of $\autp{\Gamma}$ on $E$ is generated by permutations $\alpha_{s,t}$ and by permutations $\gamma_i\in \symg{E}$ corresponding to the elements $\mathcal{H}_i$ of $\autp{\Gamma}$ satisfying
$\nu^*(\mathcal{H}_i)\in\symg{r,\mathcal P}$. In particular,
\begin{equation}\label{eq:horisontal}\gamma_i\colon (y_1,\dots,y_i,y_{i+1},\dots,y_r)\mapsto(y_1,\dots,y_iy_{i+1}y_i^{-1},y_i,\dots,y_r),\end{equation}
where $i\in \{1,2,\dots,r\}\setminus \{r_1,\dots,r_\ell\}$ and $y_i=\eta(x_i)$, $i=1,\dots,r$.

Now we are ready to describe the action of $\autp{\Gamma}$ on the set $\Epi{\Gamma,G}$.

 \begin{theorem}\label{thm:equiactions} Let $\Gamma$ be a planar Fuchsian group with presentation \eqref{eq:planarfuchsian} and let $\grp{G}$ be a finite group. For $\eta\in \Epi{\Gamma,\grp{G}}$
 set $\vec{\eta}=(\eta(x_1),\dots,\eta(x_r))$. 
 Two epimorphisms $\eta_1,\eta_2\in\Epi{\Gamma,\grp{G}}$ are equivalent if and only if there exist $a\in\aut{\grp{G}}$, $\pi\in \langle \alpha_{s,t}\rangle$, $1\leq s<t\leq r$, and $\gamma\in \langle\gamma_i\rangle$, $i\in \{1,2,\dots,r\}\setminus \{r_1,\dots,r_\ell\}$,
 such that $a(\vec\eta_1)=\gamma\pi(\vec\eta_2)$. 
\end{theorem}

\begin{proof} By Theorem~\ref{Lloyd}, $\eta_1$ is equivalent to $\eta_2$ if and only if there
exist $a\in \aut{\grp{G}}$ and $\alpha\in\autp{\Gamma}$ such that
$a(\vec\eta_1)=\vec\eta\alpha$. Since $P_r^*\lhd \autp{\Gamma}$, the group decomposes into a disjoint
union of cosets $\autp{\Gamma}=\cup_{\mu\in\symg{r,{\mathcal P}}} x_\mu P_r^*$, where $x_\mu\in (\nu^*)^{-1}(\mu)$.
Therefore, the action of $\alpha$ is determined by the action of $\bar{x}_\mu\cdot\pi$, for some $\mu\in\symg{r}$ and $\pi\in \langle \alpha_{s,t}\rangle$, $1\leq s<t\leq r$. By $\bar x_\mu$ we have denoted a transformation
of $\Epi{\Gamma,G}$ that corresponds to $x_\mu$.   
Denote by $\tau_i=(i,i+1)\in \symg{r}$ the transposition swapping the points $i$ and $i+1$.
Each $\mu\in\symg{r,{\mathcal P}}$ can be expressed as a composition $\mu=\tau_{i_1}\tau_{i_2}\dots\tau_{i_k}$,
where $\{i_1,i_2,\dots,i_k\}\subseteq  \{1,2,\dots,r\}\setminus \{r_1,\dots,r_\ell\}$. Recall that $\omega\colon \tilde {\mathcal F}_r\to \Gamma$ is the natural projection taking $\tilde x_i\mapsto x_i$, for $i=1,2,\dots,r$.
Then $x_\mu$ can be expressed as the product $\prod_{j=1}^k \mathcal{H}_{i_j}$. Setting $\bar x_\mu=\gamma=\prod_{j=1}^k\gamma_{i_j}$ we obtain the statement.   
\end{proof}

\begin{example}
Let $$\Gamma=\langle x_1,x_2,x_3\ |\ x_1^{m_1}=x_2^{m_2}=x_3^{m_3}=x_1x_2x_3=1 \rangle$$ be a triangle group with signature $(0; m_1,m_2,m_3)$. We claim that the vertical action on $\Epi{\Gamma,\grp{G}}$ coincides with the action of the group $\operatorname{Inn}(\Gamma)$ of inner automorphisms. Indeed, the vertical action is generated by the automorphisms $\mathcal{A}_{1,2}$, $\mathcal{A}_{1,3}$, and $\mathcal{A}_{2,3}$, where the images of the generators are defined by the following table:
\begin{table}[H]
\centering
\begin{tabular}{|l|lll|}
\hline\hline
& $x_1$ & $x_2$ & $x_3$ \\
\hline
\rule{0pt}{14pt}$\mathcal{A}_{1,2}$ & $x_1^{(x_1x_2)}$ & $x_2^{x_1}$ & $x_3$ \\
\rule{0pt}{14pt}$\mathcal{A}_{1,3}$ & $x_1^{(x_1x_3)}$ & $x_2^{[x_1,x_3]}$ & $x_3^{x_1}$ \\
\rule[-5pt]{0pt}{16pt}$\mathcal{A}_{2,3}$ & $x_1$ & $x_2^{(x_2x_3)}$ & $x_3^{x_2}$ \\
\hline\hline
\end{tabular}
\caption{Automorphisms of a triangle group}\label{hypermap}
\end{table}
For a generator $x_i \in \Gamma$, denote by $I_{x_i}$ the inner automorphism determined by conjugation by $x_i$. Checking the Table~\ref{hypermap}, we see that $\mathcal{A}_{1,2}=I_{x_1}I_{x_2}$, $\mathcal{A}_{1,3}=I_{x_1}I_{x_3}$, $\mathcal{A}_{2,3}=I_{x_2}I_{x_3}$. It follows that
the action of the outer automorphism group $\outp{\Gamma}$ on the set of smooth epimorphisms $\Gamma\to\grp{G}$ is isomorphic to
the action of a subgroup of $\symg{3}$. 
An epimorphic image of $\Gamma$ by an smooth epimorphism in a finite group $\grp{G}$ can be identified with a \emph{regular hypermap} $\mathcal{M}$ of type $\{m_1,m_2,m_3\}$. In particular, there is a bijection between the kernels of epimorphisms in $\Epi{\Gamma,G}$ and the isomorphism classes of hypermaps of type $\{m_1,m_2,m_3\}$. Regular hypermaps up to genus $101$ were determined by Conder \cite{Co2, Co4, Co3, Co1}. Hypermaps in the Conder's lists are determined up to isomorphism and up to geometric dualities. Two hypermaps are isomorphic if and only if the images of the three generators of $\Gamma$ are simultaneously conjugated in $\grp{G}$. In view of the previous discussion, this happens if and only if the corresponding vectors of the generators are equivalent with respect to the vertical action. The correspondence 
between the items in Conder's lists and the equivalence classes is not perfect. While Conder lists the hypermaps
up to dualities and up to taking mirror image, the corresponding actions may or may not be equivalent.
For example, there are four isomorphism classes of chiral maps of genus $7$ that split into two pairs of the form: a map and its mirror image.
The corresponding actions have the signatures $(0;2,6,9)$ and $(0;2,7,7)$. In the first case, the action of the group $\autp{\Gamma}$ is trivial over the set of isomorphism classes of the maps, and
thus the two respective maps correspond with two classes of topological equivalence. In the second case, the action of $\autp{\Gamma}$ on the isomorphism classes is non-trivial and the two chiral maps of type $\{7,7\}$, distinguished just by orientation, belong to the same equivalence class.
 
\end{example}

The following observation is well-known; see Broughton~\cite{broughton1991classifying}.

\begin{corollary}\label{cor:abelian}
If $G$ is abelian, then the vertical action is always trivial. Consequently,  $\autp{\Gamma}$ acts on $\Epi{\Gamma, \grp{G}}$ as the group of $\mathcal{P}$-invariant permutations that permute the indices
of the corresponding vectors.
\end{corollary}
\begin{proof} Since $G$ is abelian, the permutations $\alpha_{s,t}$ are trivial. Furthermore, permutations~$\gamma_i$ justs swaps the generators $y_i$ and $y_{i+1}$ of $\grp{G}$.
\end{proof}

\begin{example} In the previous text, we have seen examples where the vertical action on $\Epi{\Gamma, \grp{G}}$ has been trivial or equal to the action of the group of inner automorphisms. In what follows,  for each genus greater than $2$, we give examples of actions, where the pure braid group action is essential.

Let $\grp{D}_{2n}=\langle x,y \mid x^n=y^2=(xy)^2=1\rangle$ be a presentation of the dihedral group of order $2n$.
Let $\Gamma$ be the planar Fuchsian group with signature $(0;2,2,n,n)$. One can easily check that 
the vectors $(x^{-2}y,y,x,x)$ and $(y,y,x^{-1},x)$ determine smooth epimorphisms $\Gamma\to\grp{D}_{2n}$.
 By the Riemann existence theorem (Theorem~\ref{thm:existence}) the corresponding actions of $\grp{D}_{2n}$ exist.
Using the Riemann-Hurwitz equation, we deduce that $\grp{D}_{2n}$ acts on a surface of 
genus $n-1$. 
Clearly, no automorphism of $\grp{D}_{2n}$ takes the subsequence $(x,x)$ onto $(x^{-1},x)$;
therefore, the corresponding epimorphisms cannot be $\aut{\grp{G}}$-equivalent. Observe that the horizontal action  $\langle\gamma_1,\gamma_3\rangle$, defined in~\eqref{eq:horisontal}, preserves
the subsequence $(x,x)$. Thus the two epimorphisms are not equivalent with respect to horizontal action. On the other hand, $\alpha_{2,3}(x^{-2}y,y,x,x)=(x^{-2}y,x^{-2}y,x^{-1},x)$ and
$(x^{-2}y,x^{-2}y,x^{-1},x)^x=(y,y,x^{-1},x)$. By Theorem~\ref{thm:equiactions}, the two actions are equivalent. 
 
The corresponding discrete actions of dihedral groups of genus $\leq 8$ have identifiers O2.8, O3.22,
O4.27, O5.44, O6.48, O7.91, O8.67, see census~\cite{karaweb}.
\end{example}

\section{Automorphisms of non-planar Fuchsian groups}
\label{sec:nonplanar}
\noindent{}In what follows, we introduce several families of automorphisms of a Fuchsian group $\Gamma$.
They are defined by the images of generators of $\Gamma$ 
with presentation (\ref{eq:fuchs}), each of the missing generators in the expression of an automorphism is fixed.

\medskip
\noindent {\bf Automorphisms $\mathcal{H}_j$:}
\begin{equation}\label{eq:Hs}
\begin{aligned}
\mathcal{H}_i\colon & x_i \mapsto x_i x_{i+1} x_i^{-1},\quad x_{i+1}\mapsto x_i,\\ & \text{for } i=1,\ldots,r-1, \text{ and }|x_i|= |x_{i+1}|
\end{aligned}
\end{equation} 
These automorphisms come from the action of braid group. It is straightforward
to verify that  $\mathcal{H}_i$ is an automorphism for every admissible $i$.
They apply if $r\geq 2$.

\medskip
\noindent {\bf Automorphisms ${\mathcal B}_j$:}
\begin{equation}\label{eq:Bs} 
\begin{aligned}
\mathcal{B}_j \colon & 
a_j\mapsto a_{j+1},\quad b_j\mapsto b_{j+1},\quad a_{j+1}\mapsto d_{j+1}^{-1}a_jd_{j+1},\\
& b_{j+1}\mapsto d_{j+1}^{-1}b_jd_{j+1},\\
&\text{where } d_{j+1}=a_{j+1}b_{j+1}a^{-1}_{j+1}b^{-1}_{j+1},\\
&\text{for }j\in\{1,\ldots,\gamma-1\}.
\end{aligned}
\end{equation} 
This family of automorphisms is introduced by Harvey in~\cite[p.~396-397]{harvey1971}. They apply if $\gamma\geq 2$.
 
\medskip
\noindent {\bf Automorphisms ${\mathcal U}_j$:} 
 \begin{equation}\label{eq:Us}
\begin{aligned}
\mathcal{U}_1 \colon & a_1\mapsto a_1b_1\\
\mathcal{U}_2 \colon & a_1\mapsto a_1b_1,\quad b_1\mapsto a_1^{-1}\\
\mathcal{U}_3 \colon & x_i\mapsto a_2x_ia_2^{-1},\text{ for }i=1,\ldots,r,\\
                     & a_1\mapsto a_2a_1,\quad a_2\mapsto b_1a_2b_1^{-1},\quad b_2\mapsto a_2b_2a_2^{-1}b_1^{-1},\\ 
                     & a_j\mapsto a_2a_ja_2^{-1},\quad b_j\mapsto a_2b_ja_2^{-1}\text{ for }j=3,\ldots,\gamma \\
\mathcal{U}_4 \colon & x_r\mapsto a_1^{-1}x_ra_1,\quad a_1\mapsto a_1^{-1}x_r^{-1}a_1x_ra_1,\quad b_1\mapsto b_1a_1^{-1}x_ra_1\\
\mathcal{U}^*_4 \colon & x_1\mapsto b_\gamma^{-1}x_1b_\gamma,\quad a_\gamma\mapsto a_\gamma b_\gamma^{-1}x_1^{-1}b_\gamma,\quad b_\gamma\mapsto b_\gamma^{-1}x_1 b_\gamma x_1^{-1}b_\gamma.
\end{aligned}
\end{equation}
Harvey referred  the automorphisms $\mathcal{U}_i$, $i=1,2,3,4$, in~\cite[p.~396-397]{harvey1971}. We derived the additional automorphism ${\mathcal U}_4^*$, a twin  to $\mathcal{U}_4$.
The automorphisms $\mathcal{U}_1$ and $\mathcal{U}_2$ apply if $\gamma\geq 1$.
The automorphisms  $\mathcal{U}_3$, $\mathcal{U}_4$ and ${\mathcal U}_4^*$
require $\gamma\geq 2$. 
The group $\langle {\mathcal U}_1, {\mathcal U}_2\rangle\cong \autp{\pi(\sur{S}_1})\cong \autp{\grp{\Gamma}(1;-)}$.
 Taking into account that $\pi(\sur{S}_1)\cong
\mathbb{Z}\times\mathbb{Z}$ one can easily identify $\autp{\pi(\sur{S}_1})$
with the group of $2$-dimensional unimodular matrices with determinant $1$.
The automorphisms $\mathcal{U}_i$, and the automorphisms ${\mathcal B}_j$, $j=1,\dots,\gamma-1$, were used in~\cite[Theorem~14]{harvey1971} to prove that there is just one class of actions of a cyclic group provided  the associated Fuchsian group is a surface group. Note that the expression of $\mathcal{U}_3$  in
\cite[p.396]{harvey1971} is incorrect, however, it does not effect 
the proof of \cite[Theorem~14]{harvey1971}.

\medskip
\noindent {\bf Automorphisms ${\mathcal C}_j$:}
\begin{equation}\label{eq:Cs}
\begin{aligned}
\mathcal{C}_1 \colon & a_1\mapsto a_1b_1^{-1}\\
\mathcal{C}_2 \colon & b_1\mapsto b_1a_1\\
\mathcal{C}_3 \colon & a_1\mapsto a_1b_1^{-1}a_2b_2a_2^{-1},\quad a_2\mapsto a_2b_2^{-1}a_2b_1a_2,\\
 & b_1\mapsto a_2b_2^{-1}a_2^{-1}b_1a_2b_2a_2^{-1}\\
 \mathcal{C}_4 \colon & a_2\mapsto a_2b_2^{-1}\\
 \mathcal{C}_5 \colon & b_2\mapsto b_2a_2.
\end{aligned}
\end{equation} 
The automorphisms $\mathcal{C}_i$, $i=1,\ldots,5$ come from McCool~\cite[p.~458]{mccool1996gen}.
The automorphisms ${\mathcal C}_i\cdot\inn{\Gamma}$, $i=1,\dots,5$
generate the outer automorphism group of the fundamental group $\Gamma$ of a double torus. They form a basis of
a set of generators of outer automorphism groups of surface groups for surfaces of genera $g\geq 2$. More precisely, we have the following

\begin{theorem}\label{thm:surg} Let $\Gamma$ be the surface
group with signature $(\gamma;-)$, $\gamma>1$. Then the  automorphism group
$\autp{\Gamma}=\langle \{{\mathcal C}_1, {\mathcal C}_2, {\mathcal C}_3, \mathcal{R}\}\cup\inn{\Gamma}\rangle$, where ${\mathcal C}_i$ for $i=1,2,3$
are defined by $(\ref{eq:Cs})$ and $\mathcal{R}$ is defined by $(\ref{eq:Rs})$. 
\end{theorem} 

\begin{proof} By \cite[p.458]{mccool1996gen}
the set  $\{\mathcal{R}^j{\mathcal C}_i\mathcal{R}^{-j}\cdot\inn{\Gamma}\}$, where $i=1,2,3$, and $j=0,\dots \gamma-1$ forms a generating set of $\outp{\Gamma}$.
\end{proof}

\begin{equation}\label{eq:Rs}
\begin{aligned}
\mathcal{R}\colon & a_j\mapsto a_{j+1},\quad b_j\mapsto b_{j+1},\text{ for } j=1,\ldots,\gamma-1,\\
& a_\gamma\mapsto a_1,\quad b_\gamma\mapsto b_1.
\end{aligned}
\end{equation}

Consequently, in the case when the investigated finite group is an image  of a surface group under a smooth epimorphism, we have  complete information to derive
the classes of topological equivalence.



\section{Algorithm}

\noindent{}Assume a finite group $H=\langle A\rangle$, $A=\{a_1,a_2,\dots,a_m\}$, acts on a set $X$. The action graph ${\mathcal G}=(X;A)$
is a directed (multi)graph with set of vertices $X$ and $(x_i,x_j)$ is a (directed) edge of $\mathcal{G}$ if there
exists $a_\ell\in A$ such that $x_j=a_\ell\cdot x_i$. Every finite action graph determines
an equivalence on $X$ defined by the connectivity components. In particular, let $(\Gamma,\grp{G})$ be a $g$-admissible pair of groups. Let $X$ be the set of $\aut{\grp{G}}$-classes of smooth epimorphisms $\Gamma\to\grp{G}$, and $A\subset \autp{\Gamma}$ be a finite set of automorphisms. Denote by $[\eta]$ the $\aut{\grp{G}}$-class containing a smooth epimorphism $\eta:\Gamma\to\grp{G}$. Then we have an action of $\langle A\rangle$ defined as follows:
$\alpha[\eta_1]=[\eta_2]$, for $\alpha\in \langle A\rangle<\autp{\Gamma}$   if and only if
$[\eta_2]=[\eta_1\alpha]$. If $\langle A\cup \inn{\Gamma}\rangle=\autp{\Gamma}$
the equivalence relation defined by the action graph $(X;A)$ coincides
with the topological equivalence. The latter  holds if and only if the set 
$\{\alpha\cdot\inn{\Gamma}\ |\ \alpha\in A\}$ generates $\outp{\Gamma}$.   

\vfill
\newpage
\paragraph{\bf Algorithm 1}
\paragraph{\bf Input:} the genus, $g\geq 2$,
\vskip5pt\noindent
{\bf Output:} List of finite group actions of genus $g$ distinguished up to topological equivalence.
For signatures with $\gamma=0$ or $r=0$, the list is exact.
Otherwise, the resulting list determines a refinement of the topological equivalence.

\begin{enumerate}[{\bfseries Step 1}]
\item Numerically solve Riemann-Hurwitz equation; find all tuples of the form  \[\mathbf{rh} := \langle g,\gamma, ord, (m_1,m_2,\ldots,m_r)\rangle\]
such that
$0\leq\gamma\leq g$, $1\leq ord\leq 84(g-1)$, $r\leq 2g+2$, $m_i$ divides $ord$
for $i=1,\dots,r$ and $2g-2=ord\cdot(2\gamma-2+\sum_{i=1}^r(1-1/m_i))$;
\item For every numerical solution $\mathbf{rh}$ construct the presentation~\eqref{eq:fuchs} of the Fuchsian group $\Gamma = \grp{F}(\gamma; m_1,\ldots,m_r)$; 
\item{
\begin{enumerate}[a)]
\item If the library of small groups~\cite{smallgroups} contains groups of order $ord$, for every group $\grp{G}$ of order $ord$ construct all homomorphisms $\eta\colon \Gamma\to\grp{G}$ given by images of generators of $\Gamma$: \[\vec{\eta} := (\eta(x_1),\ldots,\eta(x_r),\eta(a_1),\ldots,\eta(a_\gamma),\eta(b_1),\ldots,\eta(b_\gamma));\]
(\texttt{hom<>} in \textsc{Magma}~\cite{magma} or \texttt{GroupHomomorphismByImages} in GAP~\cite{GAP});
\item either, construct all low-index normal subgroups of $\Gamma$ (of index $ord$); for every subgroup $N\unlhd\Gamma$ in the list set $\eta\colon\Gamma\to\Gamma/N$\\(\texttt{LowIndexNormalSubgroups} in \textsc{Magma}~\cite{magma} or \texttt{GQuotients} in GAP~\cite{GAP});
\end{enumerate}
}

\item For every homomorphism $\eta$, check whether elements of $\vec{\eta}$ satisfy the relations of $\Gamma$ and whether $|\eta(\Gamma)| = ord$. In this way we create a list $\mathcal{L}$ of representatives of the classes of the $\aut{G}$-equivalence, where
$\grp{G}=\eta(\Gamma)$ is fixed and $\Gamma$ is determined by $\mathbf{rh}$; a list $\mathcal{L}$ of tuples of the form $\langle \mathbf{rh}, \eta(\Gamma), \vec{\eta}\rangle$;
\item If $\mathcal{L}$ contains a unique element $\langle \mathbf{rh}, \eta(\Gamma), \vec{\eta}\rangle$, append the tuple to the census;
\item{ If $|\mathcal{L}|>1$, construct an action multigraph $(\mathcal{L}; A)$ on the set of representatives $\mathcal L$, where $A\subset \autp{\Gamma}$ is chosen as follows:
\begin{enumerate}[a)]
\item if $\gamma=0$, then 
\[A=\{\mathcal{A}_{s,t}\ |\ 1 \leq s < t \leq r\}\cup 
\{\mathcal{H}_i\ |\ i\in\{1,\dots,r-1\}\setminus\{r_1,\dots,r_\ell\}\},\] see (\ref{eq:purebraid}), (\ref{eq:horisontal}) and (\ref{eq:Hs});
\item if $\Gamma$ is a surface group set 
$A=\{\mathcal{C}_1,\mathcal{C}_2,\mathcal{C}_3,
\mathcal{R}\}$, see (\ref{eq:Bs}) and (\ref{eq:Rs}),
\item otherwise, set
\begin{align*}A=\{\mathcal{U}_1, \mathcal{U}_2,
\mathcal{U}_3, \mathcal{U}_4, \mathcal{U}^*_4\}&\cup \{\mathcal{ B}_j\ |\ j=1,\dots,\gamma-1\}\\& \cup 
\{\mathcal{H}_i\ |\ i\in\{1,\dots,r-1\}\setminus\{r_1,\dots,r_\ell\}\},\end{align*} see (\ref{eq:Us}), (\ref{eq:Bs}), and (\ref{eq:Hs});
\end{enumerate}
}
\item Compute the connectivity components of the multigraph $(\mathcal{L}; A)$ and for each
component store one representative.
\end{enumerate}

\newpage
\section{Completeness of the classification for $2\leq g\leq 9$}
\noindent{}For any genus $g\geq 2$,  Theorem~\ref{thm:equiactions} and Theorem~\ref{thm:surg}
prove that Algorithm~1 produces  a complete set of representatives of classes of topological equivalence of finite groups actions on $\sur{S}_g$
with planar signatures and with signatures of the form $(\gamma;-)$.
If both $r>0$ and $\gamma>0$, the algorithm determines orbits
of the action of $\autp{\grp G}\times H$ on $\Epi{\Gamma,G}$ where $H\leq\autp{\Gamma}$ is a subgroup defined in Step~6(c) of Algorithm~1.

\begin{theorem}\label{thm:complete} With the above notation, for  $2\leq g\leq 8$ and
for every $g$-admissible pair $(\Gamma,\grp{G})$ of groups, the orbits 
 of  $\autp{\grp G}\times \grp{H}$ on $\Epi{\Gamma,\grp{G}}$ determine the classes
 of the topological equivalence.
\end{theorem}
\begin{proof}
Assume that Theorem~\ref{thm:complete} does not hold. Then
there exists a $g$-admissible pair $(\Gamma,\grp{G})$, such that
$\autp{\grp G}\times \autp{\Gamma}$ has fewer orbits than $\autp{\grp G}\times \grp{H}$.
This can only happen, if $\gamma>0$, $r>0$ and the number
of orbits of $\autp{\grp G}\times \grp{H}$ on $\Epi{\Gamma,\grp{G}}$ is $>1$.
Checking the data produced by an implementation of Algorithm~1 we identify  seven such pairs $(\Gamma,\grp{G})$
for $g\leq 8$, see Table~\ref{census:mul}.
In all cases except O7.39 the group $G$ is abelian.
Employing Lemma~\ref{lem:seta}, one can easily check that   the vectors in Table~\ref{census:mul} represent pairwise different classes of the topological equivalence for abelian $G$. In case O7.39,  assume (to the contrary) that the two smooth epimorphisms represent the same equivalence class. Applying Lemma~\ref{lem:normaliser}
we get
$$N_\grp{G}(\langle \eta_1(x_1),\eta_1(x_2),\eta_1(x_3)\rangle)=N_\grp{G}(\langle xz,z^2,xz^{-1}\rangle)\cong \mathrm{C}_{2} \times \mathrm{C}_{2},$$ 
and
$$N_\grp{G}(\langle \eta_2(x_1),\eta_2(x_2),\eta_2(x_3)\rangle)=N_\grp{G}(\langle z^2,z^2,z^2\rangle)\cong \mathrm{C}_{2},$$
a contradiction.
\end{proof}

{\small
\begin{longtable}{|l|llr|}\caption{Multiple group actions with $\gamma>0$ and $r>0$, $2\leq g\leq 8$}
\label{census:mul} \\

\hline 
\endhead

\hline \multicolumn{4}{|r|}{{Continued on next page}} \\ \hline
\endfoot
 
\hline \hline
\endlastfoot 
 
\hline\hline
O5.11 & $(1; 2^4)$ &  & {\small $\grp{G}\cong\mathrm{C}_{2} \times \mathrm{C}_{2}$}\\
\hline
\rule{0pt}{11pt}Presentation & \multicolumn{3}{|l|}{$ \grp{G}=\langle x,y\ |\ x^2=y^2=[x,y]=1\rangle $}\\
\rule{0pt}{11pt}Actions & \multicolumn{3}{|l|}{$ (x, x, y, y, 1, 1)$}\\
        & \multicolumn{3}{|l|}{$(x, x, x, x, y, 1)$}\\
\hline
O5.25 & $(1; 2^2)$  &  & {\small $\grp{G}\cong\mathrm{C}_{4}\times \mathrm{C}_{2} $}\\
\hline
\rule{0pt}{11pt}Presentation & \multicolumn{3}{|l|}{$ \langle x\ |\ x^4=y^2=[x,y]=1\rangle $}\\
\rule{0pt}{11pt}Actions & \multicolumn{3}{|l|}{$ (y, y, x, 1)$}\\
        & \multicolumn{3}{|l|}{$(x^{2}, x^{2}, y, x)$}\\       
\hline
O7.9 & $(1; 3^6)$  &  & {\small $\grp{G}\cong\mathrm{C}_{3}$}\\
\hline
\rule{0pt}{11pt}Presentation & \multicolumn{3}{|l|}{$ \grp{G}=\langle x\ |\ x^3=1\rangle $}\\
\rule{0pt}{11pt}Actions & \multicolumn{3}{|l|}{$ (x, x, x, x, x, x,1,1)$}\\
        & \multicolumn{3}{|l|}{$(x^{-1}, x^{-1}, x^{-1},x,x,x,1,1)$}\\
\hline
O7.14 & $(1; 2^6)$  &  & {\small $\grp{G}\cong\mathrm{C}_{2} \times \mathrm{C}_{2}$}\\
\hline
\rule{0pt}{11pt}Presentation & \multicolumn{3}{|l|}{$ \grp{G}=\langle x,y,z\ |\ x^2=y^2=[x,y]= xyz=1\rangle $}\\
\rule{0pt}{11pt}Actions & \multicolumn{3}{|l|}{$ 
(x, x, y, y, y, y, 1, 1) $}\\
        & \multicolumn{3}{|l|}{$(x, x, y, y, z, z, 1, 1)$}\\
        & \multicolumn{3}{|l|}{$(x, x, x, x, x, x, y, 1)$}\\
\hline
O7.16 & $(1; 4^4)$  &  & {\small $\grp{G}\cong\mathrm{C}_{4}$}\\
\hline
\rule{0pt}{11pt}Presentation & \multicolumn{3}{|l|}{$ \grp{G}=\langle x\ |\ x^4=1\rangle $}\\
\rule{0pt}{11pt}Actions & \multicolumn{3}{|l|}{$ (x, x, x, x,1,1) $}\\
        & \multicolumn{3}{|l|}{$(x^{-1}, x^{-1}, x, x, 1,1)$}\\
        \hline
O7.39 & $(1; 2^3)$  &  & {\small $\grp{G}\cong\mathrm{D}_{8}$}\\
\hline
\rule{0pt}{11pt}Presentation & \multicolumn{3}{|l|}{$ \grp{G}=\langle x,z\ |\ x^2=z^4=(xz)^2=1\rangle $}\\
\rule{0pt}{11pt}Actions & \multicolumn{3}{|l|}{$ (xz,z^2, xz^{-1},x,1) $}\\
        & \multicolumn{3}{|l|}{$(z^{2}, z^{2}, z^{2}, xz^{-1} ,x)$}\\     
        \hline
O8.12 & $(1; 2^7)$  &  & {\small $\grp{G}\cong\mathrm{C}_{2}\times \mathrm{C}_{2}$}\\
\hline
\rule{0pt}{11pt}Presentation & \multicolumn{3}{|l|}{$ \grp{G}=\langle x,y,z\ |\ x^2=y^2=[x,y]= xyz=1\rangle  $}\\
\rule{0pt}{11pt}Actions & \multicolumn{3}{|l|}{$ (x,y,y,y,z,z,z,1,1) $}\\
        & \multicolumn{3}{|l|}{$(x,y,z,z,z,z,z,1,1)$}\\

\end{longtable}
} 

\begin{theorem}\label{thm:complete9} With the above notation, for  $g=9$ and
for every $g$-admissible pair $(\Gamma,\grp{G})$ of groups, except O9.109, the orbits 
 of  $\autp{\grp G}\times \grp{H}$ on $\Epi{\Gamma,\grp{G}}$ determine the classes
 of the topological equivalence. There are at most two classes of topological equivalence in case O9.109.
\end{theorem}

\begin{proof}
Assume to the contrary that Theorem~\ref{thm:complete9} does not hold. Then
there exists a $9$-admissible pair $(\Gamma,\grp{G})$, such that
$\autp{\grp G}\times \autp{\Gamma}$ has fewer orbits than $\autp{\grp G}\times \grp{H}$.
This can only happen, if $\gamma>0$, $r>0$ and the number
of orbits of $\autp{\grp G}\times \grp{H}$ on $\Epi{\Gamma,\grp{G}}$ is greater than one.
Checking the data produced by an implementation of Algorithm~1 we identified $13$ such pairs $(\Gamma,\grp{G})$, listed in Table~\ref{census:mul9}. 

In cases O9.10, O9.15, O9.18, O9.20, O9.26, O9.52, and O9.54 the group $\grp{G}$ is abelian.
Employing Lemma~\ref{lem:seta}, it can be checked that in each of the listed cases the corresponding vectors in Table~\ref{census:mul9} represent pairwise different classes of the topological equivalence for abelian $\grp{G}$.

In what follows we shall discuss the remaining cases:

\noindent\textbf{Case O9.50:}
Since $G$ is abelian, Lemma~\ref{lem:seta} distinguishes O9.50.1 from the other two actions. Now we shall distinguish O.50.2 from O.50.3. Denote by $\eta_1$ and $\eta_2$, the epimorphisms $\Gamma\to\langle x\rangle\times\langle y\rangle\cong\mathrm{C}_4\times\mathrm{C}_2$, determined by the vectors $(x^2, x^2, y, y, x, 1)$ and $(yx^2, yx^2, y, y, x, 1)$, respectively. Assume, to the contrary, that $\eta_1\sim\eta_2$. By Theorem~\ref{Lloyd} there exist $\alpha\in\operatorname{Aut}^+(\Gamma)$ and $a\in\operatorname{Aut}(\grp{G})$ such that 
$\eta_1\alpha\mid_{\{x_1,x_2,x_3,x_4\}} = a\eta_2\mid_{\{x_1,x_2,x_3,x_4\}}$. Since $\grp{G}$ is abelian, on the left-hand side we get a vector that is a re-arrangement of $(x^2,x^2,y,y)$. It follows that either $a(y)=y$ and $a(yx^2)=x^2$, or $a(y)=x^2$ and $a(yx^2)=y$. There are exactly three involutions in \grp{G}, namely $y$, $x^2$, and $yx^2$, arranged in the orbits $\{x^2\}$ and $\{y,yx^2\}$ in the action of  $\aut{\grp{G}}$. It follows that $a^{-1}(x^2)=x^2$, a contradiction with  $a^{-1}(x^2)\in\{y,yx^2\}$.

\noindent\textbf{Case O9.51:}
For each of the three actions we compute the normaliser of the subgroup generated by the images of the elliptic elements by the corresponding epimorphisms $\Gamma\to\grp{G}$. We get $N_\grp{G}(\langle x,y\rangle) \cong \grp{D}_8$, $N_\grp{G}(\langle x\rangle) \cong \grp{C}_2$, and $N_\grp{G}(\langle (xy)^2,x\rangle) \cong \grp{C}_2\times \grp{C}_2$, respectively.
By Lemma~\ref{lem:normaliser}, the actions are pairwise non-equivalent.


\noindent\textbf{Case O9.69:}
By Table~\ref{census:mul9}, there are two actions to consider.
Denote by $\eta_1$ and $\eta_2$, the respective epimorphisms $\Gamma\to\langle x, y\ |\ x^5=y^2=(xy)^2=1\rangle\cong\grp{D}_{10}$, determined by the vectors $(x, x^{-1}, y, 1)$ and $(x, x^{-2}, x^{-2}, y)$. Assume to the contrary $\eta_1\sim\eta_2$. By Theorem~\ref{Lloyd} there exist $\alpha\in\operatorname{Aut}^+(\Gamma)$ and $a\in\operatorname{Aut}(\grp{G})$ such that
$a\eta_1\alpha\mid_{\{x_1,x_2\}} = \eta_2\mid_{\{x_1,x_2\}}$. By~\cite[Theorem~8]{zieschang1966automorphismen}, every automorphism $\alpha\in\autp{\Gamma}$ takes an elliptic element to a conjugate of an elliptic element. Taking into the account $\grp{G}\cong\grp{D}_{10}$ we get:
$$(x,x^{-2})\in \{ (a(x),a(x^{-1}), (a(x^{-1}),a(x)), (a(x),a(x)), (a(x^{-1}),a(x^{-1}))\}.$$
Since $x$ is element of order five, $x\neq x^{-1}$. Thus $(x,x^{-2})\neq (a(x),a(x))$, and
$(x,x^{-2})\neq (a(x^{-1}),a(x^{-1}))$. It follows that
either $(a(x),a(x^{-1})) = (x,x^{-2})$ or $(a(x^{-1}),a(x)) = (x,x^{-2})$. Clearly, $a(x)=x^i$ for $i\in\{1,2,3,4\}$. We conclude that $(x^{\pm i},x^{\mp i})=(x,x^{-2})$, which is impossible.



\noindent\textbf{Case O9.102:}
For the two actions we compute the normaliser of the subgroup generated by the images of the elliptic elements by the corresponding epimorphisms $\Gamma\to\grp{G}$. We get $|N_\grp{G}(\langle y\rangle)| = 8$, while $|N_\grp{G}(\langle (yx)^2,x^2\rangle)| = 16$.
By Lemma~\ref{lem:normaliser}, the actions are non-equivalent.

\noindent\textbf{Case O9.104:}
Denote by $\eta_1$ and $\eta_2$, the epimorphisms $\Gamma\to\langle x\rangle\times\langle y\rangle\cong\mathrm{C}_8\times\mathrm{C}_2$, determined by the vectors $(y, y, x, 1)$ and $(x^4, x^4, y, x)$, respectively. Assume that $\eta_1\sim\eta_2$. By Theorem~\ref{Lloyd} there exist $\alpha\in\operatorname{Aut}^+(\Gamma)$ and $a\in\operatorname{Aut}(\grp{G})$ such that 
$\eta_1\alpha\mid_{\{x_1,x_2\}} = a\eta_2\mid_{\{x_1,x_2\}}$. Since $\grp{G}$ is abelian, on the left-hand side we get the vector $(y,y)$, while on the right-hand side we get $(a(x^4),a(x^4))$. It follows that $a(x^4)=y$. There are exactly three involutions in \grp{G}, namely $y$, $x^4$, and $yx^4$, arranged in the orbits $\{x^4\}$ and $\{y,yx^4\}$ in the action of  $\aut{\grp{G}}$. Similarly as above, $x^4=a^{-1}(y)$ implying $y=x^4$, a contradiction. 
\end{proof}

{\small
\begin{longtable}{|l|llr|}\caption{Multiple group actions with $\gamma>0$ and $r>0$, $ g=9$}
\label{census:mul9} \\

\hline 
\endhead

\hline \multicolumn{4}{|r|}{{Continued on next page}} \\ \hline
\endfoot
 
\hline \hline
\endlastfoot 
\hline
\hline
 O9.10 & $(1; 3^8)$ &  & {\small $\grp{G}\cong\mathrm{C}_{3}$}\\
\hline
\rule{0pt}{11pt}Presentation & \multicolumn{3}{|l|}{$ \grp{G}=\langle x\ |\ x^3=1\rangle $}\\
\rule{0pt}{11pt}Actions & \multicolumn{3}{|l|}{O9.10.1 $ (x, x^{-1}, x^{-1}, x^{-1}, x^{-1}, x^{-1}, x^{-1}, x^{-1},1,1)$}\\
        & \multicolumn{3}{|l|}{O9.10.2 $(x,x,x,x, x^{-1}, x^{-1}, x^{-1}, x^{-1},1,1)$}\\
\hline
 O9.15 & $(2; 2^4)$ &  & {\small $\grp{G}\cong\mathrm{C}_{2}\times\mathrm{C}_2$}\\
\hline
\rule{0pt}{11pt}Presentation & \multicolumn{3}{|l|}{$ \grp{G}=\langle x,y\ |\ x^2=y^2=1,[x,y]=1\rangle $}\\
\rule{0pt}{11pt}Actions & \multicolumn{3}{|l|}{O9.15.1 $ (x, x, y, y,1,1,1,1)$}\\
        & \multicolumn{3}{|l|}{O9.15.2 $(x,x,x,x,y,1,1,1)$}\\
\hline
 O9.18 & $(1; 2^8)$ &  & {\small $\grp{G}\cong\mathrm{C}_{2}\times\mathrm{C}_2$}\\
\hline
\rule{0pt}{11pt}Presentation & \multicolumn{3}{|l|}{$ \grp{G}=\langle x,y,z\ |\ x^2=y^2=z^2=1,[x,y]=1,xyz=1\rangle $}\\
\rule{0pt}{11pt}Actions & \multicolumn{3}{|l|}{O9.18.1 $(x,x,y,y,y,y,y,y,1,1)$}\\
        & \multicolumn{3}{|l|}{O9.18.2 $(x,x,x,x,y,y,y,y,1,1)$}\\
        & \multicolumn{3}{|l|}{O9.18.3 $(z,z,x,x,y,y,y,y,1,1)$}\\
        & \multicolumn{3}{|l|}{O9.18.4 $(x,x,x,x,x,x,x,x,y,1)$}\\
\hline
 O9.20 & $(1; 2^2,4^4)$ &  & {\small $\grp{G}\cong\mathrm{C}_{4}$}\\
\hline
\rule{0pt}{11pt}Presentation & \multicolumn{3}{|l|}{$ \grp{G}=\langle x\ |\ x^4=1\rangle $}\\
\rule{0pt}{11pt}Actions & \multicolumn{3}{|l|}{O9.20.1 $(x^2, x^2, x, x, x, x, 1, 1)$}\\
        & \multicolumn{3}{|l|}{O9.20.2 $(x^2, x^2, x^{-1}, x^{-1}, x, x, 1, 1)$}\\
\hline
 O9.26 & $(1; 5^4)$ &  & {\small $\grp{G}\cong\mathrm{C}_{5}$}\\
\hline
\rule{0pt}{11pt}Presentation & \multicolumn{3}{|l|}{$ \grp{G}=\langle x\ |\ x^5=1\rangle $}\\
\rule{0pt}{11pt}Actions & \multicolumn{3}{|l|}{O9.26.1 $(x^3,x^{-1},x^{-1},x^{-1},1,1)$}\\
        & \multicolumn{3}{|l|}{O9.26.2 $(x,x,x^{-1},x^{-1},1,1)$}\\
        & \multicolumn{3}{|l|}{O9.26.3 $(x^2,x^3,x,x^{-1},1,1)$}\\
\hline
 O9.50 & $(1; 2^4)$ &  & {\small $\grp{G}\cong\mathrm{C}_{4}\times\mathrm{C}_2$}\\
\hline
\rule{0pt}{11pt}Presentation & \multicolumn{3}{|l|}{$ \grp{G}=\langle x,y\ |\ x^4=y^2=1,[x,y]=1\rangle $}\\
\rule{0pt}{11pt}Actions & \multicolumn{3}{|l|}{O9.50.1 $(y, y, y, y, x, 1)$}\\
        & \multicolumn{3}{|l|}{O9.50.2 $(x^2, x^2, y, y, x, 1)$}\\
        & \multicolumn{3}{|l|}{O9.50.3 $(yx^2, yx^2, y, y, x, 1)$}\\
\hline
 O9.51 & $(1; 2^4)$ &  & {\small $\grp{G}\cong\mathrm{D}_{8}$}\\
\hline
\rule{0pt}{11pt}Presentation & \multicolumn{3}{|l|}{$ \grp{G}=\langle x,y\ |\ x^2=y^2=(xy)^4=1\rangle $}\\
\rule{0pt}{11pt}Actions & \multicolumn{3}{|l|}{O9.51.1 $ (x, x, y,y, 1,1)$}\\
        & \multicolumn{3}{|l|}{O9.51.2 $(x, x, x, x, y, 1)$}\\
        & \multicolumn{3}{|l|}{O9.51.3 $((xy)^2, (xy)^2, x, x, y, 1)$}\\
\hline
 O9.52 & $(1; 2^4)$ &  & {\small $\grp{G}\cong\mathrm{C}_{2}\times\mathrm{C}_{2}\times\mathrm{C}_2$}\\
\hline
\rule{0pt}{11pt}Presentation & \multicolumn{3}{|l|}{$ \grp{G}=\langle x,y,z\ |\ x^2=y^2=z^2=1,[x,y]=[x,z]=[y,z]=1\rangle $}\\
\rule{0pt}{11pt}Actions & \multicolumn{3}{|l|}{O9.52.1 $(xyz, x, y, z, 1, 1)$}\\
        & \multicolumn{3}{|l|}{O9.52.2 $(x, x, y, y, z, 1)$}\\
        & \multicolumn{3}{|l|}{O9.52.3 $(x, x, x, x, y, z)$}\\
\hline
 O9.54 & $(1; 2,4^2)$ &  & {\small $\grp{G}\cong\mathrm{C}_{4}\times\mathrm{C}_2$}\\
\hline
\rule{0pt}{11pt}Presentation & \multicolumn{3}{|l|}{$ \grp{G}=\langle x,y\ |\ x^4=y^2=1,[x,y]=1\rangle $}\\
\rule{0pt}{11pt}Actions & \multicolumn{3}{|l|}{O9.54.1 $ (y, x, yx^{-1}, 1,1)$}\\
        & \multicolumn{3}{|l|}{O9.54.2 $(x^2, x, x, yx^{-1},1)$}\\
\hline
 O9.69 & $(1; 5^2)$ &  & {\small $\grp{G}\cong\mathrm{D}_{10}$}\\
\hline
\rule{0pt}{11pt}Presentation & \multicolumn{3}{|l|}{$ \grp{G}=\langle x,y\ |\ x^5=y^2=(xy)^2=1\rangle $}\\
\rule{0pt}{11pt}Actions & \multicolumn{3}{|l|}{O9.69.1 $ (x, x^{-1}, y, 1)$}\\
        & \multicolumn{3}{|l|}{O9.69.2 $(x, x^{-2}, x^{-2}, y)$}\\
\hline
 O9.102 & $(1; 2^2)$ &  & {\small $\grp{G}\cong(\mathrm{C}_{4}\times\mathrm{C}_2)\rtimes\mathrm{C_2}$}\\
\hline
\rule{0pt}{11pt}Presentation & \multicolumn{3}{|l|}{$ \grp{G}=\langle x,y,z\ |\ x^4=y^2=1,[x,y]=1, (zx)^2=(zy)^2=1\rangle $}\\
\rule{0pt}{11pt}Actions & \multicolumn{3}{|l|}{O9.102.1 $ (y, y, x, 1)$}\\
        & \multicolumn{3}{|l|}{O9.102.2 $((yx)^2, x^{2}, y, x)$}\\
\hline
 O9.104 & $(1; 2^2)$ &  & {\small $\grp{G}\cong\mathrm{C}_{8}\times\mathrm{C}_2$}\\
\hline
\rule{0pt}{11pt}Presentation & \multicolumn{3}{|l|}{$ \grp{G}=\langle x,y\ |\ x^8=y^2=1,[x^2,y]=1,(xy)^4=1\rangle $}\\
\rule{0pt}{11pt}Actions & \multicolumn{3}{|l|}{O9.104.1 $ (y, y, x, 1)$}\\
        & \multicolumn{3}{|l|}{O9.104.2 $(x^4,x^4,y,x)$}\\
\hline
 O9.109 & $(1; 2^2)$ &  & {\small $\grp{G}\cong\mathrm{C}_{2}\times\mathrm{D}_8$}\\
\hline
\rule{0pt}{11pt}Presentation & \multicolumn{3}{|l|}{$ \grp{G}=\langle x,y,z,w\ |\ x^2=y^2=z^2=w^2=1,[x,y]=[y,z]=[x,z]=1, (zw)^4=1\rangle $}\\
\rule{0pt}{11pt}Actions & \multicolumn{3}{|l|}{O9.109.1 $ (x,y,z,w)$}\\
        & \multicolumn{3}{|l|}{O9.109.2 $(wx,yw,z,w)$}\\ 
\hline

\end{longtable}
} 
\newpage
\section{Classification of finite groups of small genera}\label{sec:class}
\noindent{}The list of groups of genus at most $9$ is
attached in the Appendix. Full data containing
representatives of each class of topological equivalence in JSON format are available at 
\url{https://www.savbb.sk/~karabas/science/discactions.html} 

For every genus $g$, $2\leq g\leq 9$, the rows
correspond to $g$-admissible pairs $(\Gamma,\grp{G})$.
Each row consists of the entries with the following meaning:
\begin{itemize}
\item {\bf Ref}: a unique identifier of the form  \verb+Og.integer+, it allows us to find the representatives of the
classes of topological equivalence in the database, the identifier comes from the census at webpage~\cite{karaweb}
\item {\bf Signature}: identifies the Fuchsian group $\Gamma$,
\item {\bf SMG id}: Identifier of the group $\grp{G}$ in the library of small groups~\cite{smallgroups, magma, GAP},
\item {\bf Structure}: Structure description of $\grp{G}$ (computed by GAP~\cite{GAP}),
\item {\bf epi}: $=|\Epi{\Gamma,\grp{G}}|$, the number of smooth epimorphisms
$\Gamma\to\grp{G}$
\item $\mathbf{autg}$: the number
of orbits of the action of $\aut{\grp{G}}$ on
$\Epi{\Gamma,\grp{G}}$, note that $\mathbf{autg} = 
\mathbf{epi}/|\aut{\grp{G}}|$,
\item $\mathbf{top}$: the number of classes of topological equivalence. 
\end{itemize}
The data were obtained by two independent implementations of Algorithm~1 in \textsc{Magma}~\cite{magma}.
The correctness of the outputs has been checked in two ways. First, all multiple actions of  $g$-admissible pairs for $2\leq g\leq 9$ were checked manually. Secondly, we have compared the obtained data with the published lists of group actions. Independent comparison of outputs was used to check for
cyclic actions. For cyclic groups the function enumerating the number of smooth 
epimorphisms is known, see \cite{liskovets2010multivariate, mednykh2006enumeration}.  

We shall discuss our findings below. By the phrase ``an admissible pair of groups
$(\Gamma,\grp{G})$ has a multiple action on a surface $\sur{S}_g$'' we mean that there are at least two
classes of topological equivalence of the actions of $G$ on $\sur{S}_g$ with signature determining $\Gamma$.
\vskip5pt
\paragraph{\bf g=2} There are $21$ admissible pairs $(\Gamma,\grp{G})$, each defines a unique equivalence class with respect to the topological equivalence. In particular, 
$\sim_{top}=\sim_{{\mathcal H}^1}$ for actions of genus $2$.
Our list of classes coincides with the one derived by Broughton in \cite{broughton1991classifying}.
\vskip5pt
\paragraph{\bf g=3} There are $49$ admissible pairs
$(\Gamma,\grp{G})$ giving rise to $55$ classes of topological equivalence. Four pairs give rise to two 
non-equivalent actions, namely, O3.9, O3.11, O3.16 and O.23. The pair $(\Gamma(0;2^2,4^2),\mathrm{C}_{4} \times \mathrm{C}_{2})$ gives rise to three
non-equivalent actions. In the list derived by Broughton
\cite[p.255]{broughton1991classifying} one admissible pair, namely,
$(\Gamma(0; 2,3,12),\mathrm{SL}(2,3) \colon \mathrm{C}_{2})$ is missing (O3.47 in our notation). It gives rise to one equivalence class. There is also apparent misprint in the signature of item 3.ag~\cite[p.255]{broughton1991classifying}. The correct signature should be $(0;3,3,7)$, O3.37 in our list.
In all other cases our data fits with the data derived by Broughton. Comparison with \cite[p.91]{breuer2000characters} 
shows
that $\sim_{top}=\sim_{{\mathcal H}^1}$ provided we are restricted to actions of genus $3$.
\vskip5pt
\paragraph{\bf g=4} There are $64$ admissible pairs
$(\Gamma,\grp{G})$ giving rise to $75$ classes of topological equivalence. Our list of actions fit
with the one derived by Bogopolski in~\cite{bogopolski1997classification}.
Genus four is the smallest genus, where the equivalences
$\sim_{top}$ and $\sim_{{\mathcal H}^1}$ do not coincide. Breuer in
\cite[p.91]{breuer2000characters} identified $73$ classes with respect to $\sim_{{\mathcal H}^1}$. The difference between
$\sim_{top}$ and $\sim_{{\mathcal H}^1}$ is explained by
Broughton in \cite[p.4]{broughton2022equivalence} as follows:
``The smallest genus example where the Abelian 
differentials fail to distinguish
topologically distinct actions occurs in genus 4 with $\grp{G} = \mathrm{C}_5$ and $4$ branch
points.'' In our list the corresponding item is
O4.13, where the pair $(\Gamma(0;5^4),\mathrm{C}_5)$ gives rise to three actions of genus four
that are not distinguished by $\sim_{{\mathcal H}^1}$. Kimura in \cite{Kimurag4} presented another
list of actions of genus four, distinguished up to relation 
$\sim_{{\mathcal H}^1}$. Breuer~\cite[p.90]{breuer2000characters}   noted that one class
is doubled, but another one is missing in the Kimura's list.
\vskip5pt
\paragraph{\bf g=5} There are $93$ admissible pairs
$(\Gamma,\grp{G})$ giving rise to $118$ classes of topological equivalence.
Breuer~\cite[p.90]{breuer2000characters} identifies $116$ classes of $\sim_{{\mathcal H}^1}$.
Hence some of the classes of topological equivalence are not distinguished by $\sim_{{\mathcal H}^1}$. 
The number of $\sim_{{\mathcal H}^1}$ classes fits with the one derived in \cite{Kuribayashi}.
Maximum multiplicity of actions of topological equivalence is five and it is achieved by the pair $(\Gamma(0;2^3,4^2), \mathrm{C}_{4} \times \mathrm{C}_{2})$ identified as O5.30 in our list. The largest order of a group with multiple action is $48$. In particular, a group of order $48$ (identified as O5.82), acts with signature $(0;3,4,4)$ in two different ways. In the list derived by Conder \cite{Co3} there are two non-isomorphic regular hypermaps corresponding to these actions, (RPH5.3 and RPH5.4).
\vskip5pt
\paragraph{\bf g=6} There are $87$ admissible pairs
$(\Gamma,\grp{G})$ giving rise to $108$ classes of topological equivalence.
Breuer~\cite[p.90]{breuer2000characters} identifies $105$ classes of $\sim_{{\mathcal H}^1}$.
Maximal multiplicity of the action is four, achieved by the actions of cyclic group $\mathrm{C}_{7}$ with signature 
$(0;7^4)$, referenced O6.26. The largest group with multiple action is the cyclic group 
$\mathrm{C}_{15}$ acting with signature $(0;5,15,15)$. The two corresponding regular hypermaps have
reference numbers RPH6.12 and RPH6.13 in \cite{Co3}.
\vskip5pt
\paragraph{\bf g=7} There are $148$ admissible pairs
$(\Gamma,\grp{G})$ giving rise to $210$ classes of topological equivalence.
Breuer~\cite[p.90]{breuer2000characters} identifies $208$ classes of $\sim_{{\mathcal H}^1}$.
Maximal multiplicity is eight achieved by the actions of $\mathrm{C}_{4}\times\mathrm{C}_{2}$ with signature 
$(0;2^4,4^2)$, with the identifier O7.48. The largest group with multiple action has order $54$ 
acting with signature $(0;2,6,9)$. It gives rise to two non-equivalent actions.
The corresponding regular map is chiral, therefore it appears
in two non-isomorphic enantiomers, reference C7.1 in Conder's list~\cite{Co2}.
\vskip5pt
\paragraph{\bf g=8} There are $108$ admissible pairs
$(\Gamma,\grp{G})$ giving rise to $150$ classes of topological equivalence.
Breuer~\cite[p.90]{breuer2000characters} identifies $141$ classes of $\sim_{{\mathcal H}^1}$.
Maximal multiplicity is six achieved by the actions of $\mathrm{C}_{10}$, referenced O8.44.
The largest group with multiple action is $\mathrm{PSL}(3,2) \colon \mathrm{C}_{2}$ of  order $336$  
acting with signature $(0;2,3,8)$. It gives rise to two non-equivalent actions.
The corresponding regular maps have references R8.1 and R8.2 in Conder's list \cite{Co1}.
\vskip5pt
\paragraph{$\mathbf{g>8}$} Computations of classes of topological equivalence for higher genera and verification of outputs are still running (in August 2023). For most recent results see~\cite{karaweb}. For example, we have found 
that our list of admissible pairs of groups diverges from the Breuer's list~\cite[p.91]{breuer2000characters}. In particular,
we got different numbers of admissible pairs $(\Gamma,\grp{G})$ for genera $9$, $13$, $17$. Our numbers are $270$, $454$, $746$, which do not fit with corresponding values in~\cite[p.91]{breuer2000characters}.  
For instance, comparing our data with these produced by Breuer we have discovered that the pairs $(\Gamma(1;2^2),\mathrm{C}_{4}\times\mathrm{C}_{4})$ and $(\Gamma(1;2^2),\mathrm{C}_{4}\times\mathrm{C}_{2}^2)$ are missing the list of $9$-admissible pairs in Breuer's data (personal communication). 
The same problem also appears in on-line database \cite{paulhusdb}.

Table~\ref{census:sum} summarises the above discussion, \textbf{Adm}, \textbf{H1}\footnote{this value corresponds to the number in column `Orb' in~\cite[Table~13]{breuer2000characters}} and \textbf{top}
denote the number of admissible pairs $(\Gamma,\grp{G})$, the number of $\sim_{{\mathcal H}^1}$
classes and the number of $\sim_{top}$ classes, respectively.

\begin{center}
\begin{longtable}{|c|c|c|c|}\caption{The numbers of finite group actions on Riemann surfaces of genera $2\leq g\leq 9$}
\label{census:sum} \\

\hline \multicolumn{1}{|c|}{\textbf{genus}} 
& \multicolumn{1}{c|}{\textbf{Adm}}
& \multicolumn{1}{c|}{$\mathbf{}$\textbf{H1}}
& \multicolumn{1}{c|}{$\mathbf{}$\textbf{top}}
 \\ \hline 

\endhead
           
2 & 21 & 21 & 21  \\
3 & 49 & 55 & 55  \\
4 & 64 & 73 & 75  \\
5 & 93 & 116 & 118 \\
6 & 87 & 105 & 108  \\
7 & 148 & 208 & 210  \\
8 & 108 & 141 & 150  \\
9 & 270 & 428\footnote{this number is not correct in~\cite{breuer2000characters}, see the discussion above} & 441  \\
\hline
\end{longtable}
\end{center}


\section*{Acknowledgement}
\noindent{}The three authors were supported by the grant GACR 20-15576S. The first and the second authors were supported by the grant APVV-19-0308 of the Slovak Research and Development Agency and by the grant VEGA~2/0078/20. The authors are very grateful to Robert Jajcay and to the anonymous referee(s) for their useful comments.

\bibliographystyle{abbrvnat}
\bibliography{equiactions}
\newpage
\appendix
\section{Census of actions of finite groups\\on Riemann surfaces of genera 2 to 9}

\small


\begin{center}
\begin{longtable}{|l|l|l|l|r|r|r|}\caption{Census of finite group actions on Riemann surfaces of genus $g=2$}
\label{census:genus2} \\

\hline \multicolumn{1}{|c|}{\textbf{Ref}} 
& \multicolumn{1}{c|}{\textbf{Signature}}
& \multicolumn{1}{c|}{\textbf{SMG id}}
& \multicolumn{1}{c|}{\textbf{Structure}}
& \multicolumn{1}{c|}{\textbf{epi}}
& \multicolumn{1}{c|}{\textbf{autg}}
& \multicolumn{1}{c|}{\textbf{top}}
 \\ \hline 
\endfirsthead

\multicolumn{7}{c}%
{{\bfseries \tablename\ \thetable{} -- continued from previous page}} \\
\hline \multicolumn{1}{|c|}{\textbf{Ref}} 
& \multicolumn{1}{c|}{\textbf{Signature}}
& \multicolumn{1}{c|}{\textbf{SMG id}}
& \multicolumn{1}{c|}{\textbf{Structure}}
& \multicolumn{1}{c|}{\textbf{epi}}
& \multicolumn{1}{c|}{\textbf{autg}}
& \multicolumn{1}{c|}{\textbf{top}} \\ \hline 
\endhead

\hline \multicolumn{7}{|r|}{{Continued on next page}} \\ \hline
\endfoot

\hline \hline
\endlastfoot

O2.1 & $(2; -)$ & $\langle 1, 1 \rangle$ & {\small $1$} & $1$ & $1$ & $1$ \\
O2.2 & $(1; 2, 2)$ & $\langle 2, 1 \rangle$ & {\small $\mathrm{C}_{2}$} & $4$ & $4$ & $1$ \\
O2.3 & $(0; 2^6)$ & $\langle 2, 1 \rangle$ & {\small $\mathrm{C}_{2}$} & $1$ & $1$ & $1$ \\
O2.4 & $(0; 3, 3, 3, 3)$ & $\langle 3, 1 \rangle$ & {\small $\mathrm{C}_{3}$} & $6$ & $3$ & $1$ \\
O2.5 & $(0; 2^5)$ & $\langle 4, 2 \rangle$ & {\small $\mathrm{C}_{2} \times \mathrm{C}_{2}$} & $60$ & $10$ & $1$ \\
O2.6 & $(0; 2, 2, 4, 4)$ & $\langle 4, 1 \rangle$ & {\small $\mathrm{C}_{4}$} & $2$ & $1$ & $1$ \\
O2.7 & $(0; 5, 5, 5)$ & $\langle 5, 1 \rangle$ & {\small $\mathrm{C}_{5}$} & $12$ & $3$ & $1$ \\
O2.8 & $(0; 2, 2, 3, 3)$ & $\langle 6, 1 \rangle$ & {\small $\mathrm{S}_{3}$} & $12$ & $2$ & $1$ \\
\end{longtable}
\end{center}



\begin{center}
\begin{longtable}{|l|l|l|l|r|r|r|}\caption{Census of finite group actions on Riemann surfaces of genus $g=3$}
\label{census:genus3} \\

\hline \multicolumn{1}{|c|}{\textbf{Ref}} 
& \multicolumn{1}{c|}{\textbf{Signature}}
& \multicolumn{1}{c|}{\textbf{SMG id}}
& \multicolumn{1}{c|}{\textbf{Structure}}
& \multicolumn{1}{c|}{\textbf{epi}}
& \multicolumn{1}{c|}{\textbf{autg}}
& \multicolumn{1}{c|}{\textbf{top}}
 \\ \hline 
\endfirsthead

\multicolumn{7}{c}%
{{\bfseries \tablename\ \thetable{} -- continued from previous page}} \\
\hline \multicolumn{1}{|c|}{\textbf{Ref}} 
& \multicolumn{1}{c|}{\textbf{Signature}}
& \multicolumn{1}{c|}{\textbf{SMG id}}
& \multicolumn{1}{c|}{\textbf{Structure}}
& \multicolumn{1}{c|}{\textbf{epi}}
& \multicolumn{1}{c|}{\textbf{autg}}
& \multicolumn{1}{c|}{\textbf{top}} \\ \hline 
\endhead

\hline \multicolumn{7}{|r|}{{Continued on next page}} \\ \hline
\endfoot

\hline \hline
\endlastfoot

O3.1 & $(3; -)$ & $\langle 1, 1 \rangle$ & {\small $1$} & $1$ & $1$ & $1$ \\
O3.2 & $(2; -)$ & $\langle 2, 1 \rangle$ & {\small $\mathrm{C}_{2}$} & $15$ & $15$ & $1$ \\
O3.3 & $(1; 2, 2, 2, 2)$ & $\langle 2, 1 \rangle$ & {\small $\mathrm{C}_{2}$} & $4$ & $4$ & $1$ \\
O3.4 & $(0; 2^8)$ & $\langle 2, 1 \rangle$ & {\small $\mathrm{C}_{2}$} & $1$ & $1$ & $1$ \\
O3.5 & $(1; 3, 3)$ & $\langle 3, 1 \rangle$ & {\small $\mathrm{C}_{3}$} & $18$ & $9$ & $1$ \\
O3.6 & $(0; 3^5)$ & $\langle 3, 1 \rangle$ & {\small $\mathrm{C}_{3}$} & $10$ & $5$ & $1$ \\
O3.7 & $(1; 2, 2)$ & $\langle 4, 1 \rangle$ & {\small $\mathrm{C}_{4}$} & $12$ & $6$ & $1$ \\
O3.8 & $(1; 2, 2)$ & $\langle 4, 2 \rangle$ & {\small $\mathrm{C}_{2} \times \mathrm{C}_{2}$} & $36$ & $6$ & $1$ \\
\end{longtable}
\end{center}



\begin{center}
\begin{longtable}{|l|l|l|l|r|r|r|}\caption{Census of finite group actions on Riemann surfaces of genus $g=4$}
\label{census:genus4} \\

\hline \multicolumn{1}{|c|}{\textbf{Ref}} 
& \multicolumn{1}{c|}{\textbf{Signature}}
& \multicolumn{1}{c|}{\textbf{SMG id}}
& \multicolumn{1}{c|}{\textbf{Structure}}
& \multicolumn{1}{c|}{\textbf{epi}}
& \multicolumn{1}{c|}{\textbf{autg}}
& \multicolumn{1}{c|}{\textbf{top}}
 \\ \hline 
\endfirsthead

\multicolumn{7}{c}%
{{\bfseries \tablename\ \thetable{} -- continued from previous page}} \\
\hline \multicolumn{1}{|c|}{\textbf{Ref}} 
& \multicolumn{1}{c|}{\textbf{Signature}}
& \multicolumn{1}{c|}{\textbf{SMG id}}
& \multicolumn{1}{c|}{\textbf{Structure}}
& \multicolumn{1}{c|}{\textbf{epi}}
& \multicolumn{1}{c|}{\textbf{autg}}
& \multicolumn{1}{c|}{\textbf{top}} \\ \hline 
\endhead

\hline \multicolumn{7}{|r|}{{Continued on next page}} \\ \hline
\endfoot

\hline \hline
\endlastfoot

O4.1 & $(4; -)$ & $\langle 1, 1 \rangle$ & {\small $1$} & $1$ & $1$ & $1$ \\
O4.2 & $(2; 2, 2)$ & $\langle 2, 1 \rangle$ & {\small $\mathrm{C}_{2}$} & $16$ & $16$ & $1$ \\
O4.3 & $(1; 2^6)$ & $\langle 2, 1 \rangle$ & {\small $\mathrm{C}_{2}$} & $4$ & $4$ & $1$ \\
O4.4 & $(0; 2^{10})$ & $\langle 2, 1 \rangle$ & {\small $\mathrm{C}_{2}$} & $1$ & $1$ & $1$ \\
O4.5 & $(2; -)$ & $\langle 3, 1 \rangle$ & {\small $\mathrm{C}_{3}$} & $80$ & $40$ & $1$ \\
O4.6 & $(1; 3, 3, 3)$ & $\langle 3, 1 \rangle$ & {\small $\mathrm{C}_{3}$} & $18$ & $9$ & $1$ \\
O4.7 & $(0; 3^6)$ & $\langle 3, 1 \rangle$ & {\small $\mathrm{C}_{3}$} & $22$ & $11$ & $2$ \\
O4.8 & $(1; 2, 2, 2)$ & $\langle 4, 2 \rangle$ & {\small $\mathrm{C}_{2} \times \mathrm{C}_{2}$} & $96$ & $16$ & $1$ \\
\end{longtable}
\end{center}



\begin{center}
\begin{longtable}{|l|l|l|l|r|r|r|}\caption{Census of finite group actions on Riemann surfaces of genus $g=5$}
\label{census:genus5} \\

\hline \multicolumn{1}{|c|}{\textbf{Ref}} 
& \multicolumn{1}{c|}{\textbf{Signature}}
& \multicolumn{1}{c|}{\textbf{SMG id}}
& \multicolumn{1}{c|}{\textbf{Structure}}
& \multicolumn{1}{c|}{\textbf{epi}}
& \multicolumn{1}{c|}{\textbf{autg}}
& \multicolumn{1}{c|}{\textbf{top}}
 \\ \hline 
\endfirsthead

\multicolumn{7}{c}%
{{\bfseries \tablename\ \thetable{} -- continued from previous page}} \\
\hline \multicolumn{1}{|c|}{\textbf{Ref}} 
& \multicolumn{1}{c|}{\textbf{Signature}}
& \multicolumn{1}{c|}{\textbf{SMG id}}
& \multicolumn{1}{c|}{\textbf{Structure}}
& \multicolumn{1}{c|}{\textbf{epi}}
& \multicolumn{1}{c|}{\textbf{autg}}
& \multicolumn{1}{c|}{\textbf{top}} \\ \hline 
\endhead

\hline \multicolumn{7}{|r|}{{Continued on next page}} \\ \hline
\endfoot

\hline \hline
\endlastfoot

O5.1 & $(5; -)$ & $\langle 1, 1 \rangle$ & {\small $1$} & $1$ & $1$ & $1$ \\
O5.2 & $(3; -)$ & $\langle 2, 1 \rangle$ & {\small $\mathrm{C}_{2}$} & $63$ & $63$ & $1$ \\
O5.3 & $(2; 2, 2, 2, 2)$ & $\langle 2, 1 \rangle$ & {\small $\mathrm{C}_{2}$} & $16$ & $16$ & $1$ \\
O5.4 & $(1; 2^8)$ & $\langle 2, 1 \rangle$ & {\small $\mathrm{C}_{2}$} & $4$ & $4$ & $1$ \\
O5.5 & $(0; 2^{12})$ & $\langle 2, 1 \rangle$ & {\small $\mathrm{C}_{2}$} & $1$ & $1$ & $1$ \\
O5.6 & $(1; 3, 3, 3, 3)$ & $\langle 3, 1 \rangle$ & {\small $\mathrm{C}_{3}$} & $54$ & $27$ & $1$ \\
O5.7 & $(0; 3^7)$ & $\langle 3, 1 \rangle$ & {\small $\mathrm{C}_{3}$} & $42$ & $21$ & $1$ \\
O5.8 & $(2; -)$ & $\langle 4, 1 \rangle$ & {\small $\mathrm{C}_{4}$} & $240$ & $120$ & $1$ \\
\end{longtable}
\end{center}



\begin{center}
\begin{longtable}{|l|l|l|l|r|r|r|}\caption{Census of finite group actions on Riemann surfaces of genus $g=6$}
\label{census:genus6} \\

\hline \multicolumn{1}{|c|}{\textbf{Ref}} 
& \multicolumn{1}{c|}{\textbf{Signature}}
& \multicolumn{1}{c|}{\textbf{SMG id}}
& \multicolumn{1}{c|}{\textbf{Structure}}
& \multicolumn{1}{c|}{\textbf{epi}}
& \multicolumn{1}{c|}{\textbf{autg}}
& \multicolumn{1}{c|}{\textbf{top}}
 \\ \hline 
\endfirsthead

\multicolumn{7}{c}%
{{\bfseries \tablename\ \thetable{} -- continued from previous page}} \\
\hline \multicolumn{1}{|c|}{\textbf{Ref}} 
& \multicolumn{1}{c|}{\textbf{Signature}}
& \multicolumn{1}{c|}{\textbf{SMG id}}
& \multicolumn{1}{c|}{\textbf{Structure}}
& \multicolumn{1}{c|}{\textbf{epi}}
& \multicolumn{1}{c|}{\textbf{autg}}
& \multicolumn{1}{c|}{\textbf{top}} \\ \hline 
\endhead

\hline \multicolumn{7}{|r|}{{Continued on next page}} \\ \hline
\endfoot

\hline \hline
\endlastfoot

O6.1 & $(6; -)$ & $\langle 1, 1 \rangle$ & {\small $1$} & $1$ & $1$ & $1$ \\
O6.2 & $(3; 2, 2)$ & $\langle 2, 1 \rangle$ & {\small $\mathrm{C}_{2}$} & $64$ & $64$ & $1$ \\
O6.3 & $(2; 2^6)$ & $\langle 2, 1 \rangle$ & {\small $\mathrm{C}_{2}$} & $16$ & $16$ & $1$ \\
O6.4 & $(1; 2^{10})$ & $\langle 2, 1 \rangle$ & {\small $\mathrm{C}_{2}$} & $4$ & $4$ & $1$ \\
O6.5 & $(0; 2^{14})$ & $\langle 2, 1 \rangle$ & {\small $\mathrm{C}_{2}$} & $1$ & $1$ & $1$ \\
O6.6 & $(2; 3, 3)$ & $\langle 3, 1 \rangle$ & {\small $\mathrm{C}_{3}$} & $162$ & $81$ & $1$ \\
O6.7 & $(1; 3^5)$ & $\langle 3, 1 \rangle$ & {\small $\mathrm{C}_{3}$} & $90$ & $45$ & $1$ \\
O6.8 & $(0; 3^8)$ & $\langle 3, 1 \rangle$ & {\small $\mathrm{C}_{3}$} & $86$ & $43$ & $2$ \\
\end{longtable}
\end{center}



\begin{center}
\begin{longtable}{|l|l|l|l|r|r|r|}\caption{Census of finite group actions on Riemann surfaces of genus $g=7$}
\label{census:genus7} \\

\hline \multicolumn{1}{|c|}{\textbf{Ref}} 
& \multicolumn{1}{c|}{\textbf{Signature}}
& \multicolumn{1}{c|}{\textbf{SMG id}}
& \multicolumn{1}{c|}{\textbf{Structure}}
& \multicolumn{1}{c|}{\textbf{epi}}
& \multicolumn{1}{c|}{\textbf{autg}}
& \multicolumn{1}{c|}{\textbf{top}}
 \\ \hline 
\endfirsthead

\multicolumn{7}{c}%
{{\bfseries \tablename\ \thetable{} -- continued from previous page}} \\
\hline \multicolumn{1}{|c|}{\textbf{Ref}} 
& \multicolumn{1}{c|}{\textbf{Signature}}
& \multicolumn{1}{c|}{\textbf{SMG id}}
& \multicolumn{1}{c|}{\textbf{Structure}}
& \multicolumn{1}{c|}{\textbf{epi}}
& \multicolumn{1}{c|}{\textbf{autg}}
& \multicolumn{1}{c|}{\textbf{top}} \\ \hline 
\endhead

\hline \multicolumn{7}{|r|}{{Continued on next page}} \\ \hline
\endfoot

\hline \hline
\endlastfoot

O7.1 & $(7; -)$ & $\langle 1, 1 \rangle$ & {\small $1$} & $1$ & $1$ & $1$ \\
O7.2 & $(4; -)$ & $\langle 2, 1 \rangle$ & {\small $\mathrm{C}_{2}$} & $255$ & $255$ & $1$ \\
O7.3 & $(3; 2, 2, 2, 2)$ & $\langle 2, 1 \rangle$ & {\small $\mathrm{C}_{2}$} & $64$ & $64$ & $1$ \\
O7.4 & $(2; 2^8)$ & $\langle 2, 1 \rangle$ & {\small $\mathrm{C}_{2}$} & $16$ & $16$ & $1$ \\
O7.5 & $(1; 2^{12})$ & $\langle 2, 1 \rangle$ & {\small $\mathrm{C}_{2}$} & $4$ & $4$ & $1$ \\
O7.6 & $(0; 2^{16})$ & $\langle 2, 1 \rangle$ & {\small $\mathrm{C}_{2}$} & $1$ & $1$ & $1$ \\
O7.7 & $(3; -)$ & $\langle 3, 1 \rangle$ & {\small $\mathrm{C}_{3}$} & $728$ & $364$ & $1$ \\
O7.8 & $(2; 3, 3, 3)$ & $\langle 3, 1 \rangle$ & {\small $\mathrm{C}_{3}$} & $162$ & $81$ & $1$ \\
\end{longtable}
\end{center}



\begin{center}
\begin{longtable}{|l|l|l|l|r|r|r|}\caption{Census of finite group actions on Riemann surfaces of genus $g=8$}
\label{census:genus8} \\

\hline \multicolumn{1}{|c|}{\textbf{Ref}} 
& \multicolumn{1}{c|}{\textbf{Signature}}
& \multicolumn{1}{c|}{\textbf{SMG id}}
& \multicolumn{1}{c|}{\textbf{Structure}}
& \multicolumn{1}{c|}{\textbf{epi}}
& \multicolumn{1}{c|}{\textbf{autg}}
& \multicolumn{1}{c|}{\textbf{top}}
 \\ \hline 
\endfirsthead

\multicolumn{7}{c}%
{{\bfseries \tablename\ \thetable{} -- continued from previous page}} \\
\hline \multicolumn{1}{|c|}{\textbf{Ref}} 
& \multicolumn{1}{c|}{\textbf{Signature}}
& \multicolumn{1}{c|}{\textbf{SMG id}}
& \multicolumn{1}{c|}{\textbf{Structure}}
& \multicolumn{1}{c|}{\textbf{epi}}
& \multicolumn{1}{c|}{\textbf{autg}}
& \multicolumn{1}{c|}{\textbf{top}} \\ \hline 
\endhead

\hline \multicolumn{7}{|r|}{{Continued on next page}} \\ \hline
\endfoot

\hline \hline
\endlastfoot

O8.1 & $(8; -)$ & $\langle 1, 1 \rangle$ & {\small $1$} & $1$ & $1$ & $1$ \\
O8.2 & $(4; 2, 2)$ & $\langle 2, 1 \rangle$ & {\small $\mathrm{C}_{2}$} & $256$ & $256$ & $1$ \\
O8.3 & $(3; 2^6)$ & $\langle 2, 1 \rangle$ & {\small $\mathrm{C}_{2}$} & $64$ & $64$ & $1$ \\
O8.4 & $(2; 2^{10})$ & $\langle 2, 1 \rangle$ & {\small $\mathrm{C}_{2}$} & $16$ & $16$ & $1$ \\
O8.5 & $(1; 2^{14})$ & $\langle 2, 1 \rangle$ & {\small $\mathrm{C}_{2}$} & $4$ & $4$ & $1$ \\
O8.6 & $(0; 2^{18})$ & $\langle 2, 1 \rangle$ & {\small $\mathrm{C}_{2}$} & $1$ & $1$ & $1$ \\
O8.7 & $(2; 3, 3, 3, 3)$ & $\langle 3, 1 \rangle$ & {\small $\mathrm{C}_{3}$} & $486$ & $243$ & $1$ \\
O8.8 & $(1; 3^7)$ & $\langle 3, 1 \rangle$ & {\small $\mathrm{C}_{3}$} & $378$ & $189$ & $1$ \\
\end{longtable}
\end{center}



\begin{center}
\begin{longtable}{|l|l|l|l|r|r|r|}\caption{Census of finite group actions on Riemann surfaces of genus $g=9$}
\label{census:genus9} \\

\hline \multicolumn{1}{|c|}{\textbf{Ref}} 
& \multicolumn{1}{c|}{\textbf{Signature}}
& \multicolumn{1}{c|}{\textbf{SMG id}}
& \multicolumn{1}{c|}{\textbf{Structure}}
& \multicolumn{1}{c|}{\textbf{epi}}
& \multicolumn{1}{c|}{\textbf{autg}}
& \multicolumn{1}{c|}{\textbf{top}}
 \\ \hline 
\endfirsthead

\multicolumn{7}{c}%
{{\bfseries \tablename\ \thetable{} -- continued from previous page}} \\
\hline \multicolumn{1}{|c|}{\textbf{Ref}} 
& \multicolumn{1}{c|}{\textbf{Signature}}
& \multicolumn{1}{c|}{\textbf{SMG id}}
& \multicolumn{1}{c|}{\textbf{Structure}}
& \multicolumn{1}{c|}{\textbf{epi}}
& \multicolumn{1}{c|}{\textbf{autg}}
& \multicolumn{1}{c|}{\textbf{top}} \\ \hline 
\endhead

\hline \multicolumn{7}{|r|}{{Continued on next page}} \\ \hline
\endfoot

\hline \hline
\endlastfoot

O9.1 & $(9; -)$ & $\langle 1, 1 \rangle$ & {\small $1$} & $1$ & $1$ & $1$ \\
O9.2 & $(5; -)$ & $\langle 2, 1 \rangle$ & {\small $\mathrm{C}_{2}$} & $1023$ & $1023$ & $1$ \\
O9.3 & $(4; 2, 2, 2, 2)$ & $\langle 2, 1 \rangle$ & {\small $\mathrm{C}_{2}$} & $256$ & $256$ & $1$ \\
O9.4 & $(3; 2^8)$ & $\langle 2, 1 \rangle$ & {\small $\mathrm{C}_{2}$} & $64$ & $64$ & $1$ \\
O9.5 & $(2; 2^12)$ & $\langle 2, 1 \rangle$ & {\small $\mathrm{C}_{2}$} & $16$ & $16$ & $1$ \\
O9.6 & $(1; 2^16)$ & $\langle 2, 1 \rangle$ & {\small $\mathrm{C}_{2}$} & $4$ & $4$ & $1$ \\
O9.7 & $(0; 2^20)$ & $\langle 2, 1 \rangle$ & {\small $\mathrm{C}_{2}$} & $1$ & $1$ & $1$ \\
O9.8 & $(3; 3, 3)$ & $\langle 3, 1 \rangle$ & {\small $\mathrm{C}_{3}$} & $1458$ & $729$ & $1$ \\
\end{longtable}
\end{center}

\end{document}